\documentclass[12pt]{article}

\usepackage{amsfonts}
\usepackage{amssymb}
\usepackage{latexsym}
\usepackage{amsmath}
\usepackage{bibmods}

\usepackage{graphicx}
\usepackage{color}
\usepackage{bibnames}

\usepackage{amssymb}
\usepackage{amscd} 
\usepackage{color}

 \newtheorem{theorem}{Theorem}[section]
 
 \newtheorem{lemma}{Lemma}[section]
 \newtheorem{corollary}{Corollary}[section]

 \newtheorem{definition}{Definition}[section]
 \newtheorem{remark}{Remark}[section]
 \numberwithin{equation}{section}

 \newcommand{\diag}{\mathrm{diag}\,}
\newcommand{\supp}{\mathrm{supp}\,}

 \newcommand{\rbx}{\hfill{\rule{1ex}{1ex}}}

\def\tov{{\overline{t}}}
 \def\ao{{\overline{\zeta}}}
\def\o{\omega(\zeta)}
\def\ot{\omega(t)}

  \newcommand{\ro}[2]{#1=0,1,\ldots, #2}

 \newcommand{\vp}{\varphi}
\newcommand{\D}{\displaystyle}

\newcommand{\ig}{\int_\Gamma}

\newcommand{\im}{\mathrm{im}\,}

\newcommand{\re}{\textrm{Re}\,}

\newcommand{\cA}{\mathcal{A}}
\newcommand{\cB}{\mathcal{B}}

\newcommand{\cJ}{\mathcal{J}}
\newcommand{\cK}{\mathcal{K}}
\newcommand{\cL}{\mathcal{L}}
\newcommand{\cM}{\mathcal{M}}

\newcommand{\sC}{{\mathbb C}}

\newcommand{\sN}{{\mathbb N}}

\newcommand{\sR}{{\mathbb R}}

\newcommand{\sZ}{{\mathbb Z}}

\title{Spline
Galerkin methods for the Sherman-Lauricella equation on contours
with corners}

\author{Victor D. Didenko\thanks{Faculty of Science, Universiti
    Brunei Darussalam, Bandar Seri Begawan, BE1410 Brunei ({\tt
      diviol@gmail.com}), supported by the Universiti Brunei
    Darussalam under Grant  UBD/GSR/S\&T/19.}, \and
 Tao Tang\thanks{Hong Kong Baptist University,
 Faculty of Science, SCT715, Cha Chi Ming Science Tower,
  Kowloon Tong, Hong Kong ({\tt ttang@hkbu.edu.hk}), partially supported by Hong Kong
Research Grants Council CERG Grants, National Science Foundation of
China, and Hong Kong Baptist University FRG Grants.}, \and  Anh My
Vu\thanks{Faculty of Science, Universiti Brunei Darussalam,  Bandar
Seri Begawan, BE1410 Brunei ({\tt  anhmy7284@gmail.com}), supported
by the Universiti Brunei
    Darussalam under Grant  UBD/GSR/S\&T/19.} }

\begin{document}

\maketitle

\begin{abstract}
Spline Galerkin approximation methods for the
Sherman-Lau\-ri\-cel\-la integral equation on simple closed
piecewise smooth contours are studied, and necessary and sufficient
conditions for their stability are obtained. It is shown that the
method under consideration is stable if and only if certain
operators associated with the corner points of the contour are
invertible. Numerical experiments demonstrate a good convergence of
the spline Galerkin methods and validate theoretical results.
Moreover, it is shown that if all corners of the contour have
opening angles located in interval $(0.1\pi, 1.9\pi)$, then the
corresponding Galerkin method based on splines of order $0$, $1$ and
$2$ is always stable. These results are in strong contrast with the
behaviour of the Nystr\"om method, which has a number of instability
angles in the interval mentioned.
\end{abstract}



\pagestyle{myheadings} \thispagestyle{plain} \markboth{VICTOR D.
DIDENKO, TAO TANG AND ANH MY VU}{GALERKIN METHODS FOR
SHERMAN-LAURICELLA EQUATION}

  \section{Introduction}

Let $D$ be a simply connected planar domain bounded by a piecewise
smooth curve $\Gamma$. It is well known that the solution of various
boundary value problems for the biharmonic equation
\begin{equation*}
\Delta^2 u(x,y) \equiv \frac{\partial^4 u}{\partial
x^4}+2\frac{\partial^4 u}{\partial x^2\partial y^2}
+\frac{\partial^4 u}{\partial y^4} =0, \quad (x,y)\in D,
\end{equation*}
where $\Delta$ is the Laplace operator, can be constructed via
solutions of boundary integral equations. Consider the biharmonic
Dirichlet problem
\begin{equation}\label{eq2}
 \begin{aligned}
 &   \Delta^2|_D=0,\\
  &  u|_\Gamma =f_1,\quad \frac{\partial u}{\partial
    \mathbf{n}}\Big|_{\Gamma} =f_2,
\end{aligned}
\end{equation}
where $\partial u/\partial \mathbf{n}$  denotes the normal
derivatives and $f_1,f_2$ are sufficiently smooth functions defined
on the boundary $\Gamma$. Setting $z=x+iy, i^2=-1$, one can identify
$D$ with a domain in the complex plane $\sC$.  This problem arises
in various applications, in particular while considering the
behaviour of viscous flows with small Reynolds numbers, bacteria
movement, deflection of plates, elastic equilibrium of solids,
sintering \cite{Crowdy2002, Greengard1996, Ka:1975, Mikh:1964,
Muskhelishvili1966, Ockendon1995, PP:1982}.

Let us equip the curve $\Gamma$ with the counterclockwise
orientation and consider the Sherman--Lauricella equation
  \begin{equation}\label{eq3}
   \omega(t)+\frac{1}{2\pi i}\int_\Gamma \omega(\zeta)\,
d\ln \left ( \frac{\zeta-t}{\overline{\zeta}-\overline{t}} \right )
- \frac{1}{2\pi i}\int_\Gamma \overline{\omega(\zeta)}\, d\, \left (
\frac{\zeta-t}{\overline{\zeta}-\overline{t}} \right ) =f(t),\quad
t=x+iy\in \Gamma,
\end{equation}
where the bar denotes the complex conjugation and $\omega$ is an
unknown function. Equation \eqref{eq3} originated in works of G.
Lauricella (see \cite{Lau:1909}). He was the first who used the
method of integral equations in elasticity.  Later D.I. Sherman
rewrites Lauricella equation in a complex form and proposes a new
simple way to derive it \cite{She:1940a}. The equation \eqref{eq3}
is uniquely solvable in appropriate functional spaces, provided $f$
satisfies certain smoothness conditions and
\begin{equation}\label{eq4}
    \re \int_\Gamma \overline{f(t)} \, dt =0,
\end{equation}
\cite{Du:1986, Mikh:1964, PP:1982}. Moreover, let
$\alpha=\alpha(x,y)$, $(x,y)\in \Gamma$ denote the angle between the
real axis $\sR$ and the outward normal $\mathbf{n}$ to $\Gamma$ at
the point $(x,y)$ and let $\mathbf{l}$ be the unit vector such that
the angle between $\mathbf{l}$ and the real axis is $\alpha-\pi/2$.
If one defines the function $f=f(t)=f(x,y)$, $t=x+iy$ by
 \begin{equation}\label{eq5}
    f(t):=e^{-i\alpha} \left ( f_2(t)+i \frac{\partial f_1}{\partial \mathbf{l}}(t)\right
    ), \quad t\in \Gamma,
\end{equation}
then the solution of the Sherman-Lauricella equation \eqref{eq3}
with such right-hand side $f$\/ can be used to determine a solution
of the boundary value problem \eqref{eq2}. More precisely, if
$\omega$ is a solution of the equation \eqref{eq3} with the
right-hand side \eqref{eq5}, consider two holomorphic functions
$\vp=\vp(z)$ and $\psi=\psi(z)$, $z\in D$ defined by
  \begin{align}
 \label{eq6}
 & \vp(z)= \frac{1}{2\pi i}\int_\Gamma\frac{\omega(\zeta)}{\zeta-z}
 \,d\zeta, \quad z\in D, \\[1ex]
  \label{eq7}
 & \psi(z)=  \frac{1}{2\pi i}\int_\Gamma\frac{\overline{\omega(\zeta)}}{\zeta-z}
 \,d\zeta + \frac{1}{2\pi i}\int_\Gamma\frac{\omega(\zeta)}{\zeta-z}
 \,d\overline{\zeta} -\frac{1}{2\pi i}\int_\Gamma\frac{\overline{\zeta}\omega(\zeta)}{(\zeta-z)^2}
 \,d\zeta, \quad z\in D.
\end{align}
According to \cite{Mikh:1964}, the pair $\{\vp(z), \psi(z)\}$
represents a solution of the boundary value problem
  \begin{equation*}
\vp(t)+t\overline{\vp'(t)} + \overline{\psi(t)} =e^{-i\alpha} \left
( f_2(t)+i \frac{\partial f_1}{\partial \mathbf{l}}(t)\right
    ),\quad t\in
\Gamma.
\end{equation*}
Therefore, by \cite[Lemma 5.1.4]{DS:2008} the function
 \begin{equation}\label{eq9}
   u(x,y):=\re(\overline{z}\vp(z)+\psi(z)), \quad z=x+iy\in D
\end{equation}
is the solution of the boundary value problem \eqref{eq2}.

Thus if an exact or an approximate solution of the integral equation
\eqref{eq3} is known, a solution of the biharmonic problem
\eqref{eq2} can be obtained by using formulas \eqref{eq6},
\eqref{eq7} and \eqref{eq9}. Therefore, the main effort should be
directed to the determination of solutions of the Sherman-Lauricella
equation \eqref{eq3}. Note that the Nystr\"om method for the
Sherman-Lauricella equation on smooth contours has been used in
\cite{Greengard1996, Kro:2001} to find approximate solution of
biharmonic problems arising in fluid dynamics. However, the authors
of those works have not presented any stability conditions for the
method considered. If $\Gamma$ has corner points, the stability
study becomes more involved since the integral operators in
\eqref{eq3} are not compact. For piecewise smooth contours,
conditions of the stability of the Nystr\"om method are established
in \cite{DH:2011b, DH:2011}. These results have been used in
\cite{DH:2013e} in order to construct a very accurate numerical
method to find solutions of the biharmonic problem \eqref{eq2} in
piecewise smooth domains in the case of piecewise continuous
boundary conditions.

In the present paper, we consider spline based Galerkin methods for
the equation \eqref{eq3} and study their stability. It is shown that
the corresponding method is stable if and only if certain operators
$R^\tau$ from an algebra of Toeplitz operators are invertible. These
operators depend on the spline space used and on the opening angles
of the corner points $\tau\in \Gamma$. Unfortunately, nowadays there
is no analytic tool to verify whether the operators in question are
invertible or not. Nevertheless, we propose a numerical approach
which can handle this problem. Thus spline Galerkin methods are
applied to the Sherman--Lauricella equation on simple model curves
and the behaviour of the corresponding approximation operators
provide an information about the invertibility of the operators
$R^\tau$, $\tau\in \Gamma$. Note that in comparison to the Nystr\"om
method, the implementation of spline Galerkin methods to solve the
Sherman--Lauricella equation, requires more preparatory work. On the
other hand, numerical experiments suggest that these methods have no
"critical" angles located in the interval $[0.1\pi, 1.9\pi]$, i.e.
if the boundary $\Gamma$ does not possess corners with opening
angles from the interval mentioned, then these methods are stable.
In a sense, this is similar to the behaviour of the corresponding
approximation methods for Sherman--Lauricella and Muskhelishvili
equations in the case of smooth curves which always converge
\cite{DV:2007, DH:2011,DS:2002}. Of course, one also has to study
the opening angles in the intervals $(0, 0.1\pi)$ and to
$(1.9\pi,2\pi)$ but this is a time consuming operation and will be
considered elsewhere.

 \section{Splines and Galerkin method}

We start this section with the construction of spline spaces on the
contour $\Gamma$. Let $\gamma=\gamma(s), s\in \sR$ be a $1$-periodic
parametrization  of $\Gamma$, and let $\cM_\Gamma$ denote the set of
all corner points $\tau_0,\tau_1, \ldots,\tau_{q-1}$ of $\Gamma$.
Without loss of generality we can assume that $\tau_j=\gamma(j/q)$
for all $j=0,1,\ldots, q-1$. In addition, we also suppose that the
function $\gamma$ is two times continuously differentiable on each
interval $(j/q,(j+1)/q)$ and
 \begin{equation*}
    \left |\gamma'\left (\frac{j}{q}+0\right) \right |=
       \left |\gamma'\left (\frac{j}{q}-0\right) \right | , \quad
       \ro j{q-1}.
\end{equation*}
Note that the last condition is not very restrictive and can always
be satisfied by changing the parametrization of $\Gamma$ in an
appropriate way.

Let  $f$ and $g$ be functions defined on the real line $\sR$, and
let $f\ast g$ denote the convolution
 $$
(f\ast g)(s):=\int_\sR f(s-x) g(x) dx
 $$
of $f$ and $g$. If $\chi$ is the characteristic function of the
interval $[0,1)$,
 $$
\chi(s):=\left \{
 \begin{array}{ll}
   1 & \text{ if } s\in [0,1), \\
   0 & \text{ otherwise}, \\
 \end{array}
\right.
 $$
then $\widehat{\phi}=\widehat{\phi}^{(d)}(s)$ refers to the function
defined by
 $$
\widehat{\phi}^{(d)}(s):=\left \{
 \begin{array}{ll}
 \chi(s) &  \text{ if } d=0 ,\\
   (\chi \ast \widehat{\phi}^{(d-1)})(s) & \text{ if } d=1,2 \ldots
   \, .
 \end{array}
 \right .
 $$
Recall that for any given non-negative integer $d$, the function
$\widehat{\phi}$ generates spline spaces on $\sR$. Thus if an
$n\in\sN$ is fixed, then closure in the $L^2$-norm of the set of all
finite linear combinations of the functions
$\widehat{\phi}_{nj}(s):=\widehat{\phi}(ns-j)$, $j\in\sZ$
constitutes a spline space on $\sR$.

Using the above defined spline functions, one can introduce spline
spaces on the contour $\Gamma$. More precisely, for a fixed
non-negative integer $d$ and an $n\in \sN$, $n\geq d+1$, we denote
by $S_n^d=S_n^d(\Gamma)$ the set of all linear combinations of the
functions
 $$
\widehat{\phi}_{nj}(t):=\widehat{\phi}(ns-j), \quad
t=\gamma(s)\in\Gamma, \quad j=0, 1, \ldots, n-(d+1), \quad s\in \sR,
 $$
the support of which belongs entirely to one of the arcs $[\tau_k,
\tau_{k+1})$, $k=0, \ldots, q$ and $\tau_{q+1}:=\tau_0$. This
definition is correct since the support
$\mathrm{supp}\,\widehat{\phi}$ of the function $\widehat{\phi}$ is
contained in the interval $[0, d+1]$ \cite{Schumaker1993} and
$\gamma$ is a $1$-periodic function.

In what follows, we also consider operators acting on various
subspaces of the Hilbert space $\widetilde{l}^2=l^2(\sZ)$ of all
sequences $(\xi_k)$ of complex numbers $\xi_k, k\in \sZ$ satisfying
the condition
 $$
||(\xi_k)||:=\left (\sum_{k\in\sZ} |\xi_k|^2 \right )^{1/2}<\infty.
 $$
The space $\widetilde{l}^2$ is closely connected to spline spaces on
the real line $\sR$. Thus the following result is true.
 \begin{lemma}[{\cite{DeBoor1978}}]\label{l1}
 Let $n\in\sN$. Then there are constants $c_1$ and $c_2$ such that
 for any sequence $(\xi_k)\in \widetilde{l}^2$ the relations
  \begin{equation*}
 ||(\xi_k)||\leq c_1 \sqrt{n}\Big|\Big|\sum_{k\in \sZ} \xi_k
 \widehat{\phi}_{nk}\Big| \Big|_{L^2(\sR)} \leq  \frac{c_2}{\sqrt{n}}||(\xi_k)||
\end{equation*}
hold.
\end{lemma}

Further, let $L^2(\Gamma)$ denote set of all Lebesgue measurable
functions $f$ such that
  $$
||f||_{L^2} :=\left ( \int_\Gamma|f(t)|^2 \,ds  \right )^{1/2}
<\infty,
  $$
and let $A_\Gamma:L^2(\Gamma)\to L^2(\Gamma)$ be the operator
corresponding to the Sherman-Lauricella equation \eqref{eq3}. It is
well known that the operator $A_\Gamma$ is not invertible on the
space $L^2(\Gamma)$ \cite{Mikh:1964}. On the other hand, the
invertibility of the corresponding operator is a necessary condition
for the applicability of any Galerkin method to any operator
equation. Therefore, for approximate solution of the equation
\eqref{eq3} we use the equation with the operator $B_{\Gamma}$
instead of $A_{\Gamma}$ and chose the right-hand sides $f$ of the
initial equation \eqref{eq3} from a suitable subspace of
$L^2(\Gamma)$. More precisely, let $W^{1,2}(\Gamma)$ denote the
closure of the set of all functions $f$ with bounded derivatives in
the norm
    $$
||f||_{W^{1,2}} :=\left ( \int_\Gamma|f(t)|^2 \,ds +
\int_\Gamma|f'(t)|^2 \,ds \right )^{1/2},
  $$
and let $T_\Gamma : L^2(\Gamma)\to L^2(\Gamma)$ refer to the
operator defined by
 \begin{equation}\label{eq10}
 T_\Gamma\omega(t):= \frac{1}{(\tov -\overline{a})}\,\frac{1}{2\pi i}\int_\Gamma\left (
\frac{\o}{(\zeta-a)^2}d\zeta +
\frac{\overline{\o}}{(\overline{\zeta}-\overline{a})^2}d\overline{\zeta}
\right ),
\end{equation}
where $a$ is a point in $D$.

 \begin{theorem}[{\cite{DH:2011}}]\label{t1}
If $\Gamma$ is a simple closed piecewise smooth contour, then
operator
\begin{equation*}
 B_\Gamma:=A_\Gamma+T_\Gamma
\end{equation*}
is invertible on the space $L_2(\Gamma)$. Moreover, if function
$f\in W_2^1(\Gamma)$ satisfies the condition \eqref{eq4}, then the
solution of the equation
 \begin{equation}\label{eq12}
B_\Gamma\omega=f
\end{equation}
belongs to the space $W_2^1(\Gamma)$ and is a solution of the
original Sherman-Lauricella equation \eqref{eq3}.
 \end{theorem}

\begin{figure}[!b]
\centering
\includegraphics[height=45mm,width=60mm]{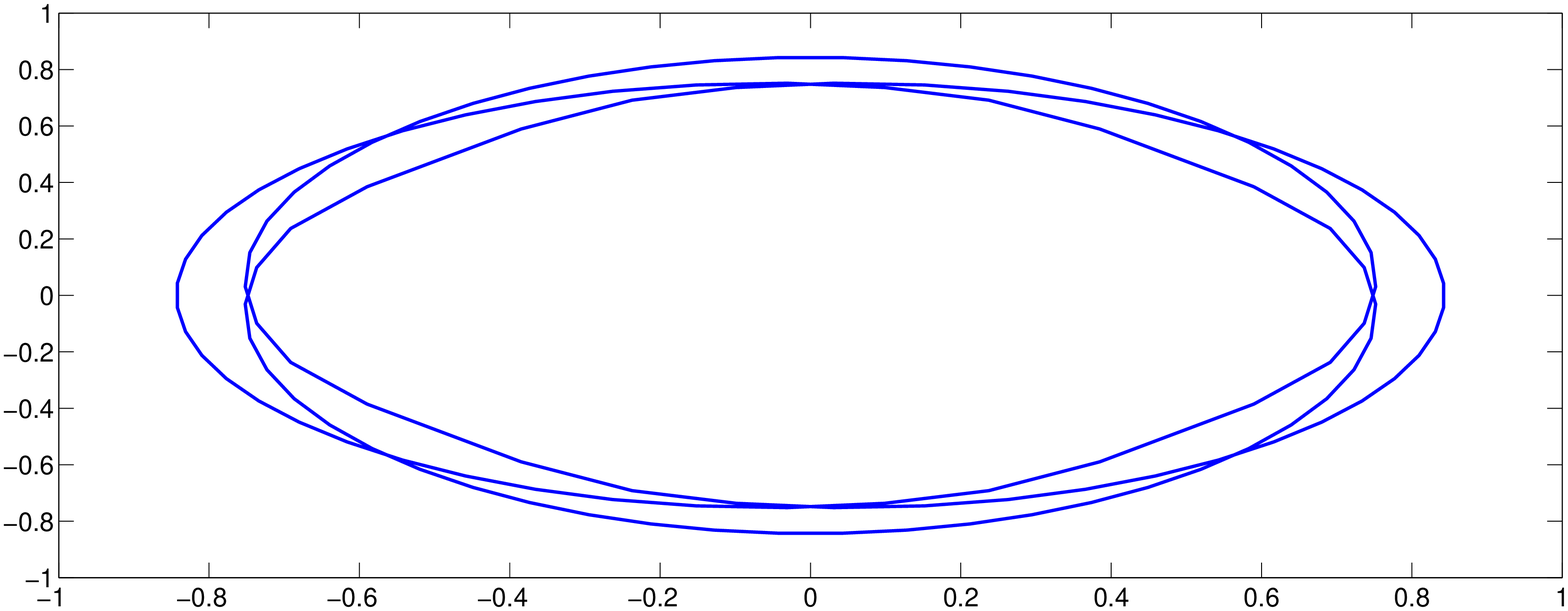}
  \includegraphics[height=45mm,width=60mm]{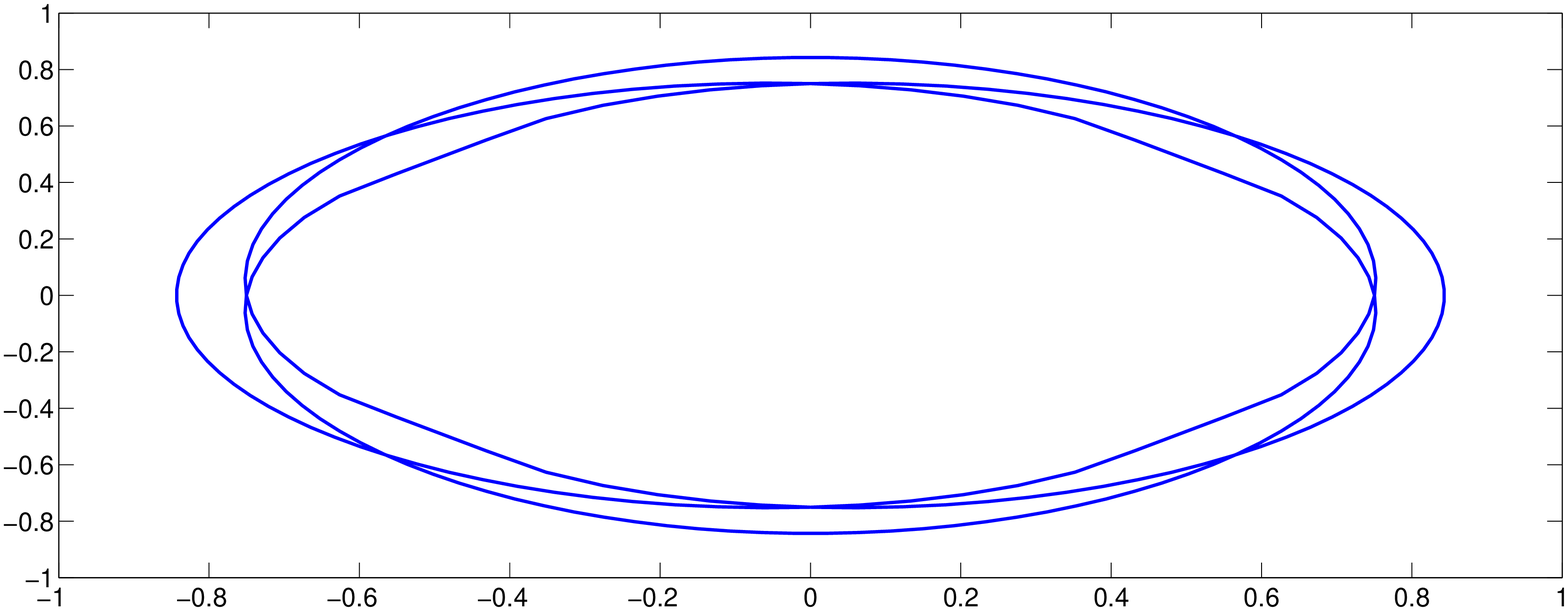}
\includegraphics[height=45mm,width=60mm]{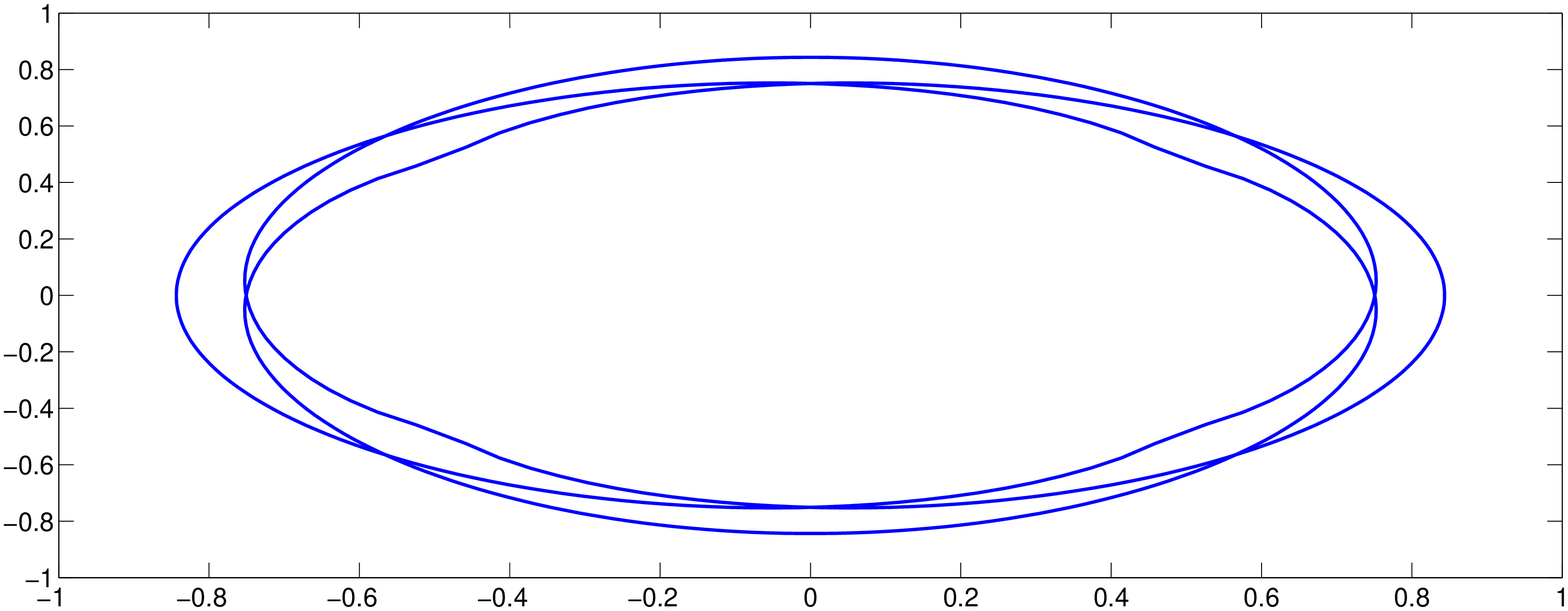}
  \includegraphics[height=45mm,width=60mm]{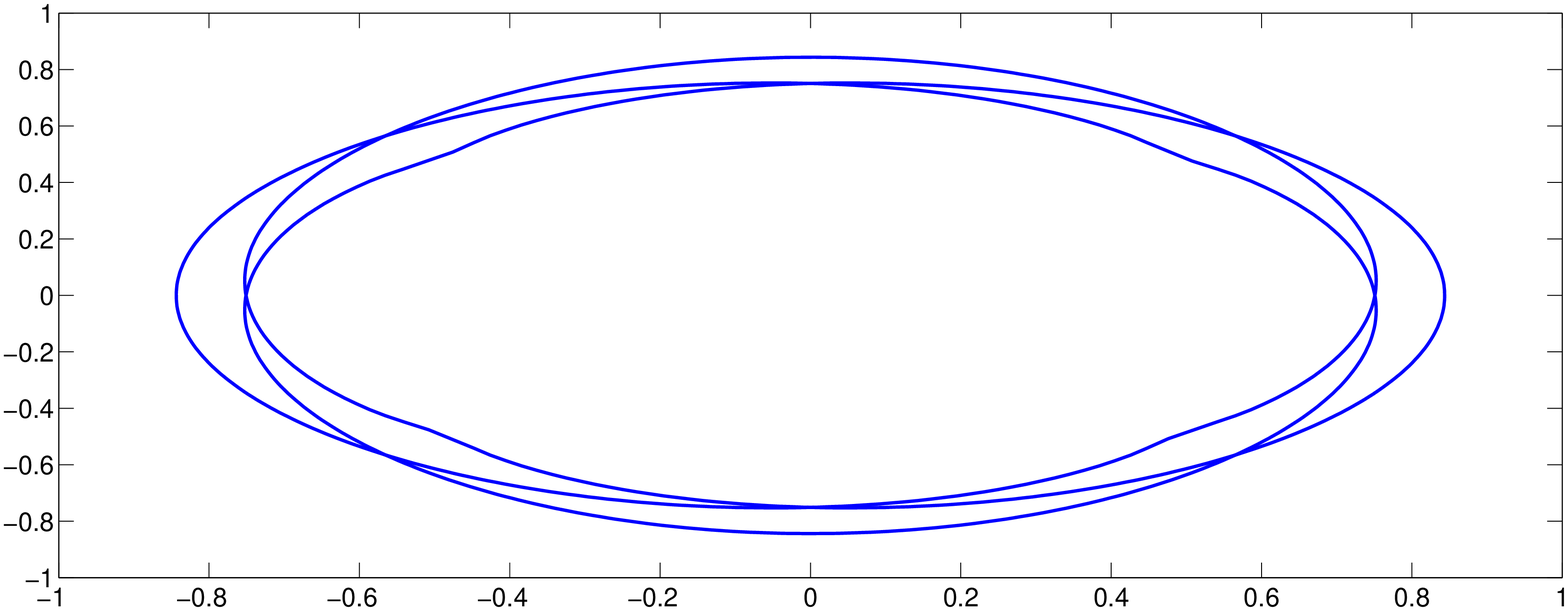}
\caption{\sf Approximate solution $\omega_n(t)$ of the
Sherman--Lauricella equation \eqref{eq3} on the unit square $\Gamma$
with $f:=f_1$ defined by \eqref{RHS} and $d=0$. From the left to the
right: $n=128,256,512,1024$} \label{fig:square}
\end{figure}

\begin{figure}[!hp]
\centering
\includegraphics[height=45mm,width=60mm]{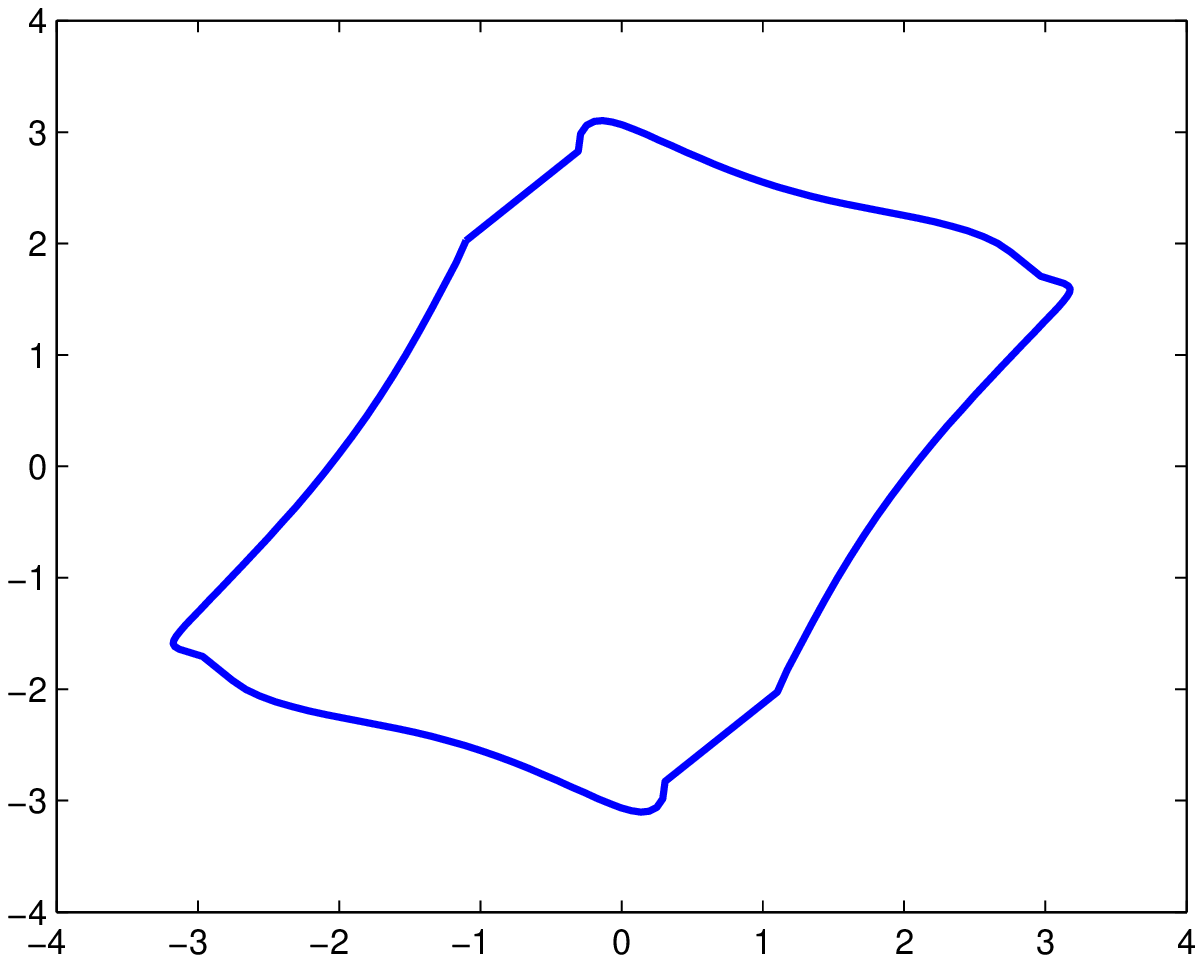}
  \includegraphics[height=45mm,width=60mm]{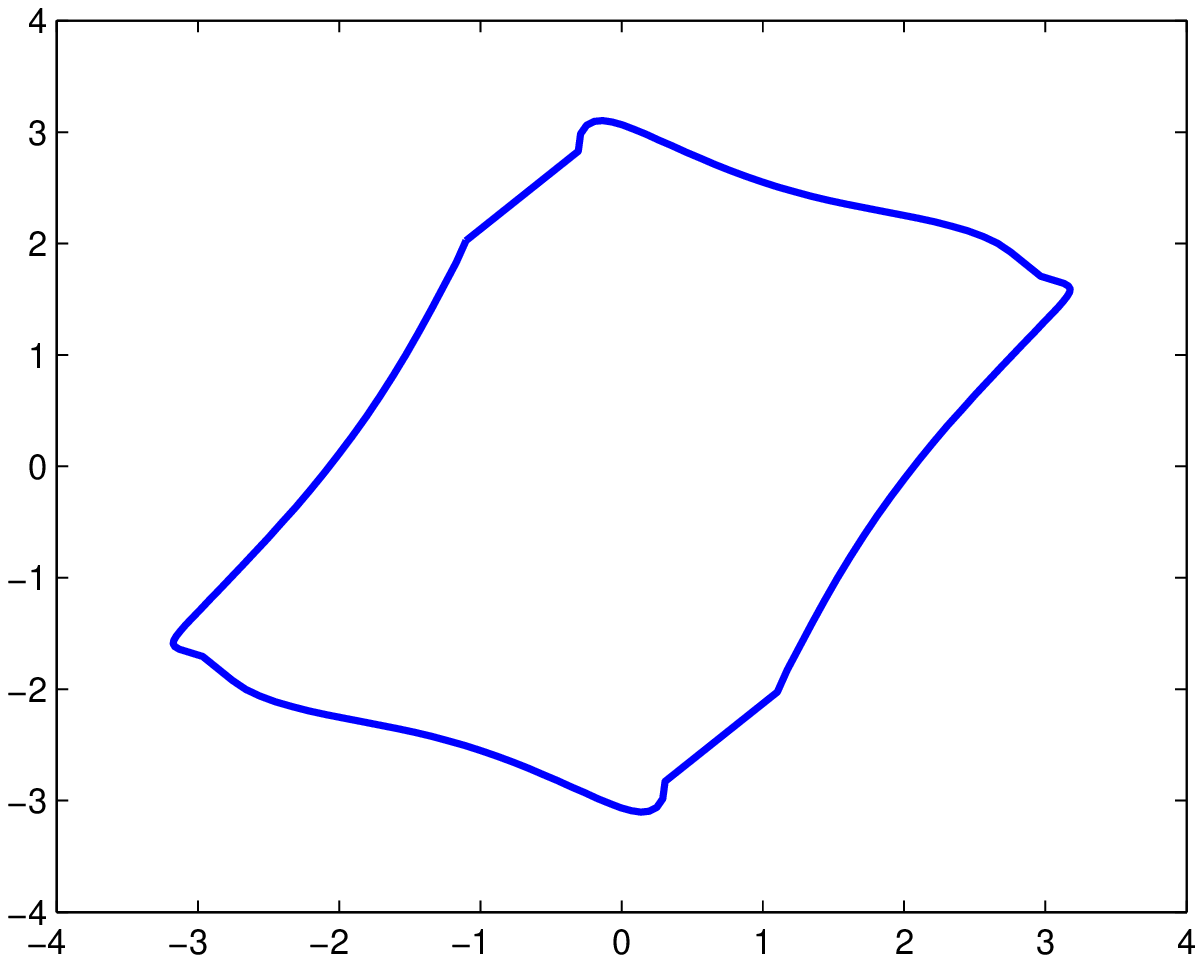}
\includegraphics[height=45mm,width=60mm]{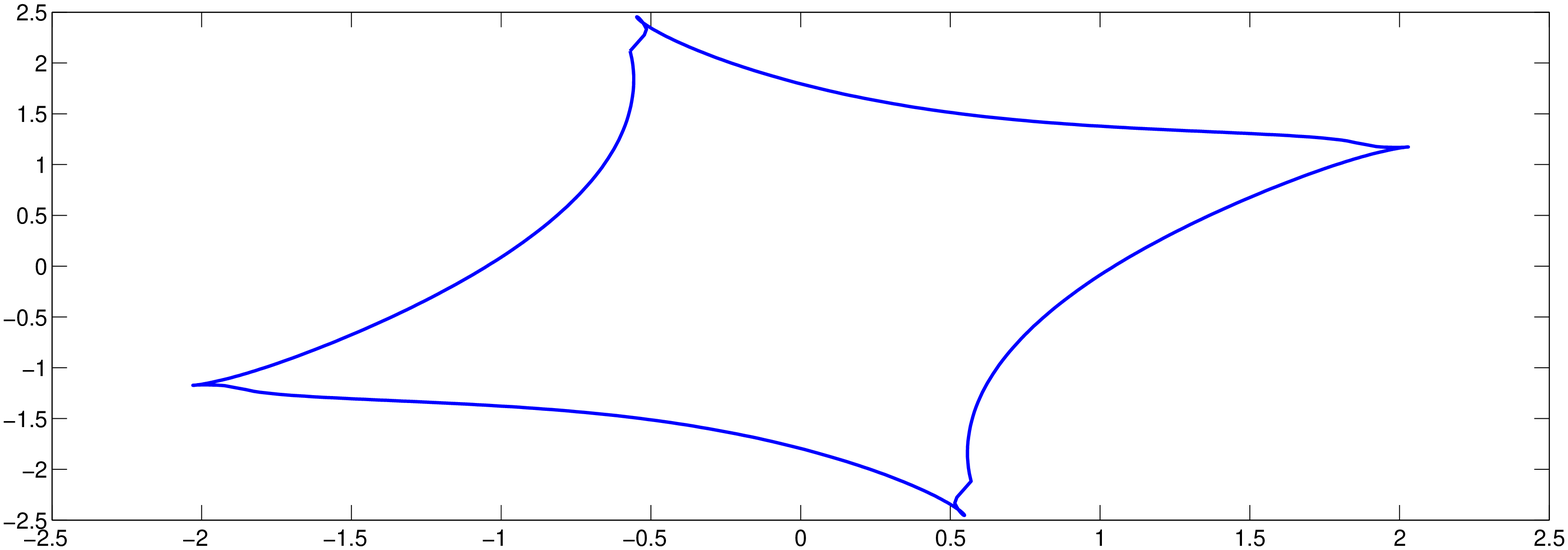}
  \includegraphics[height=45mm,width=60mm]{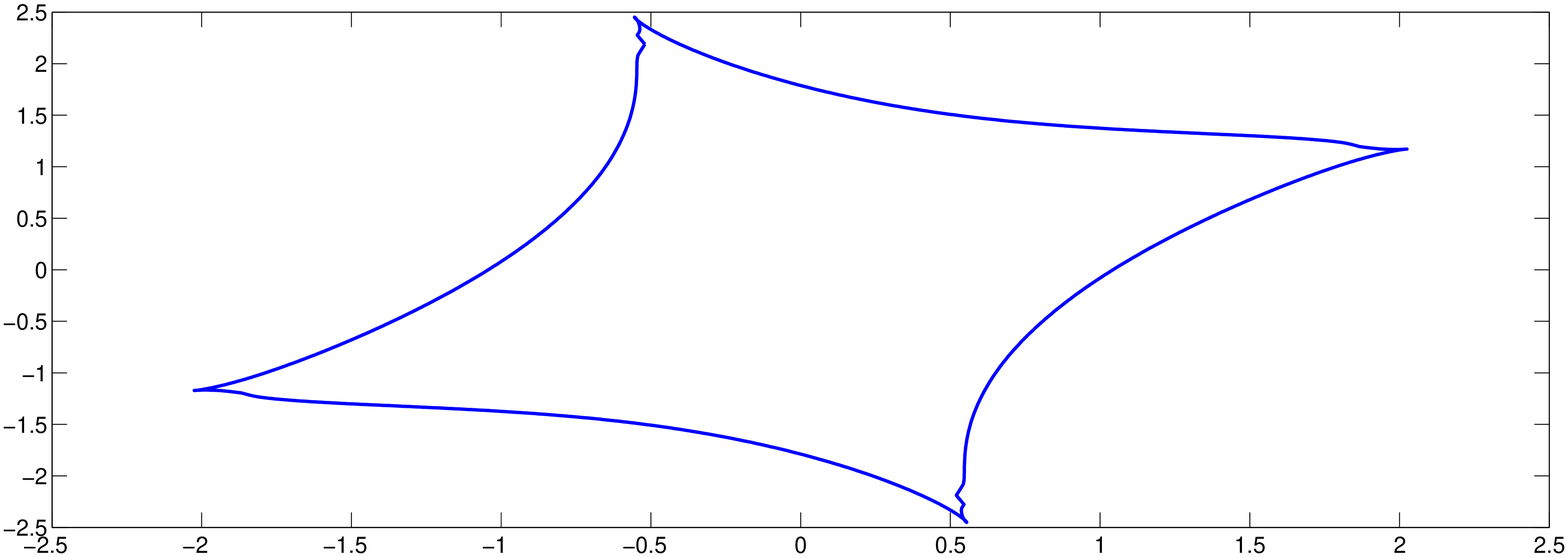}
\caption{\sf Approximate solution $\omega_n(t)$ of the
Sherman--Lauricella equation \eqref{eq3} on the rhombus $\Gamma$,
$\alpha=\pi/3$ with $f:=f_1$ defined by \eqref{RHS} and $d=0$. From
the left to the right: $n=128,256,512,1024$}
\label{fig:om_rh_unit_pi3}
\includegraphics[height=45mm,width=60mm]{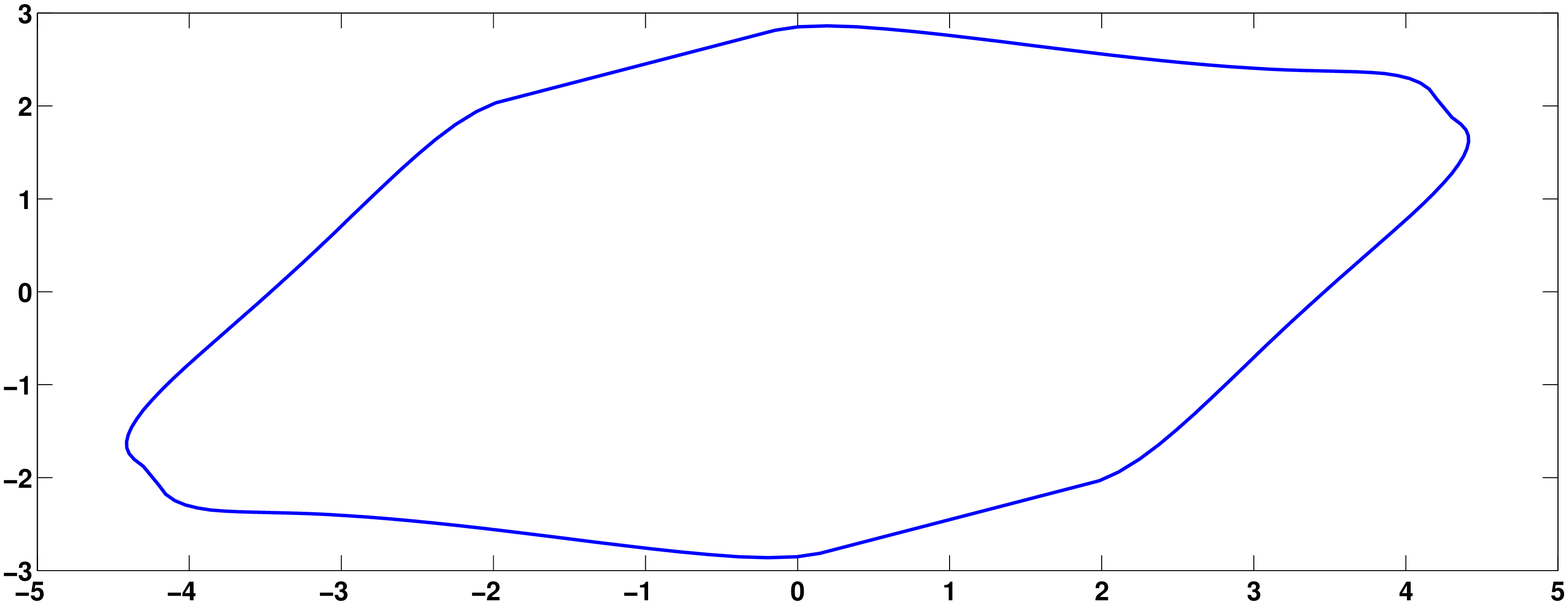}
  \includegraphics[height=45mm,width=60mm]{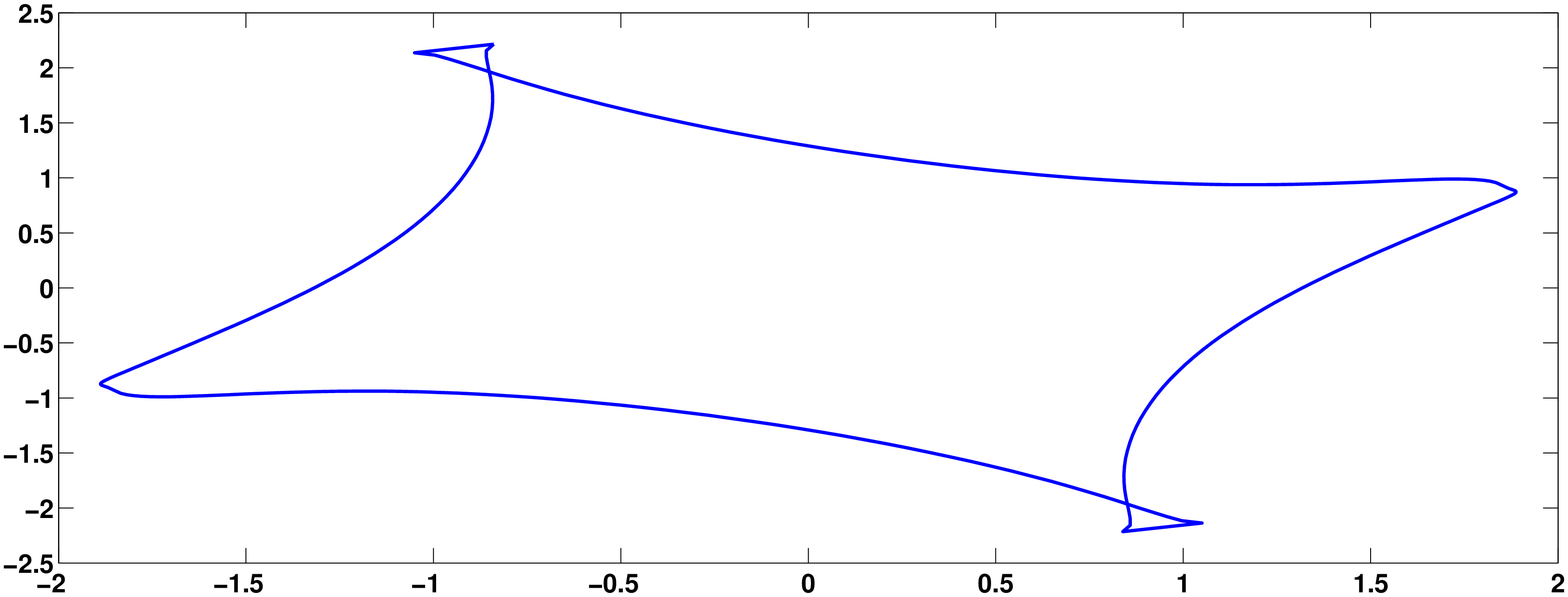}
\includegraphics[height=45mm,width=60mm]{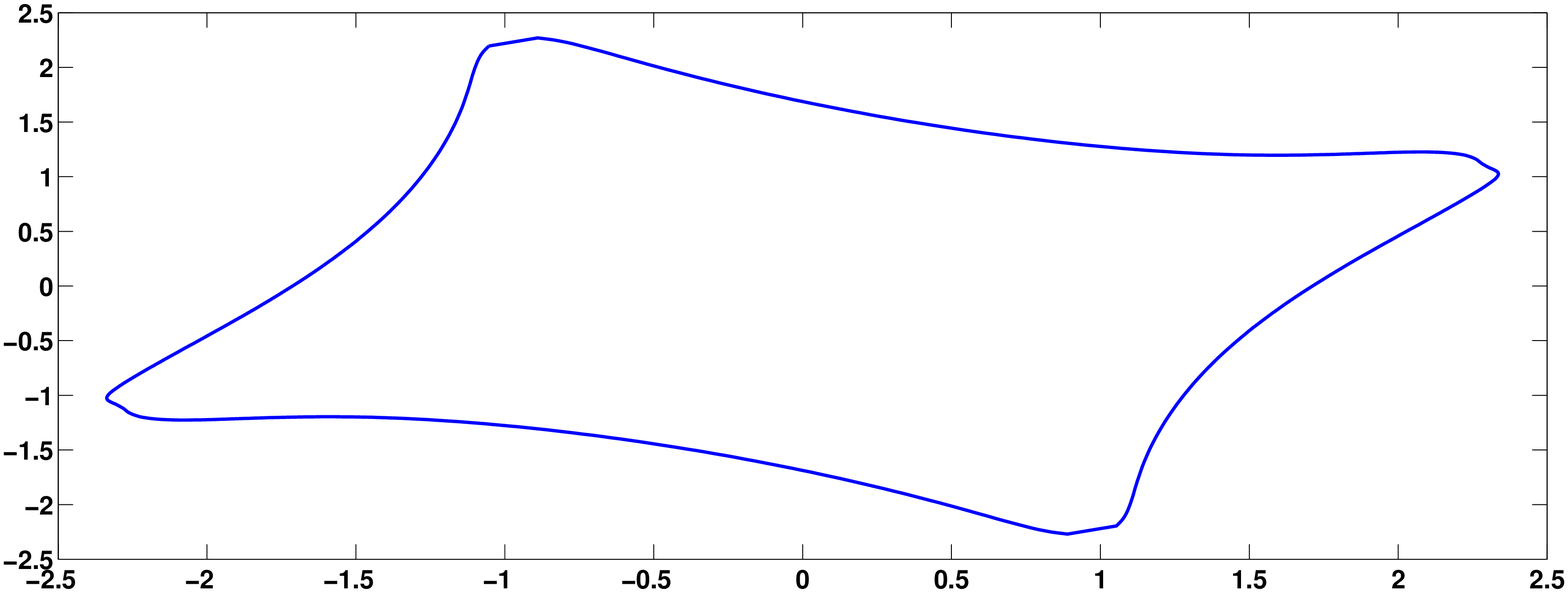}
  \includegraphics[height=45mm,width=60mm]{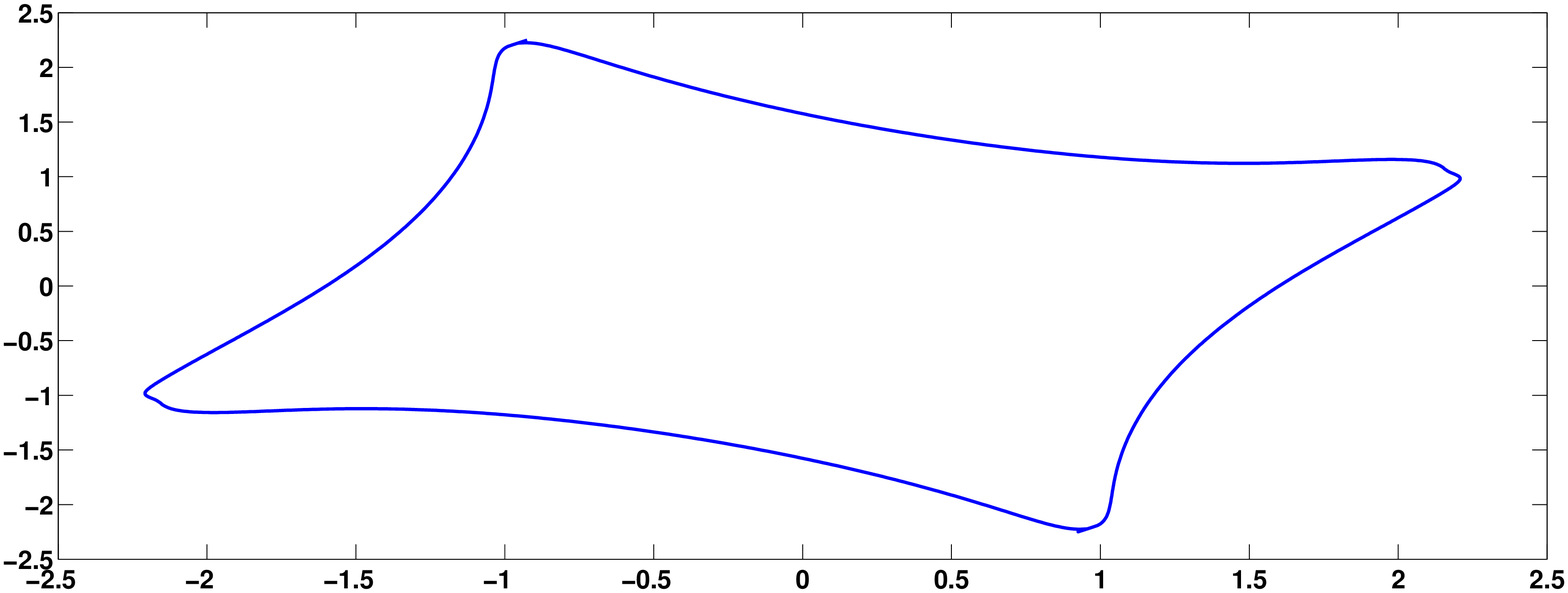}
\caption{\sf Approximate solution $\omega_n(t)$ of the
Sherman--Lauricella equation \eqref{eq3} on the rhombus $\Gamma$,
$\alpha=\pi/4$ with $f:=f_1$ defined by \eqref{RHS} and $d=0$. From
the left to the right: $n=128,256,512,1024$} \label{fig:om_rh_pi4}
\end{figure}

\begin{figure}[!ht]
\centering
\includegraphics[height=45mm,width=60mm]{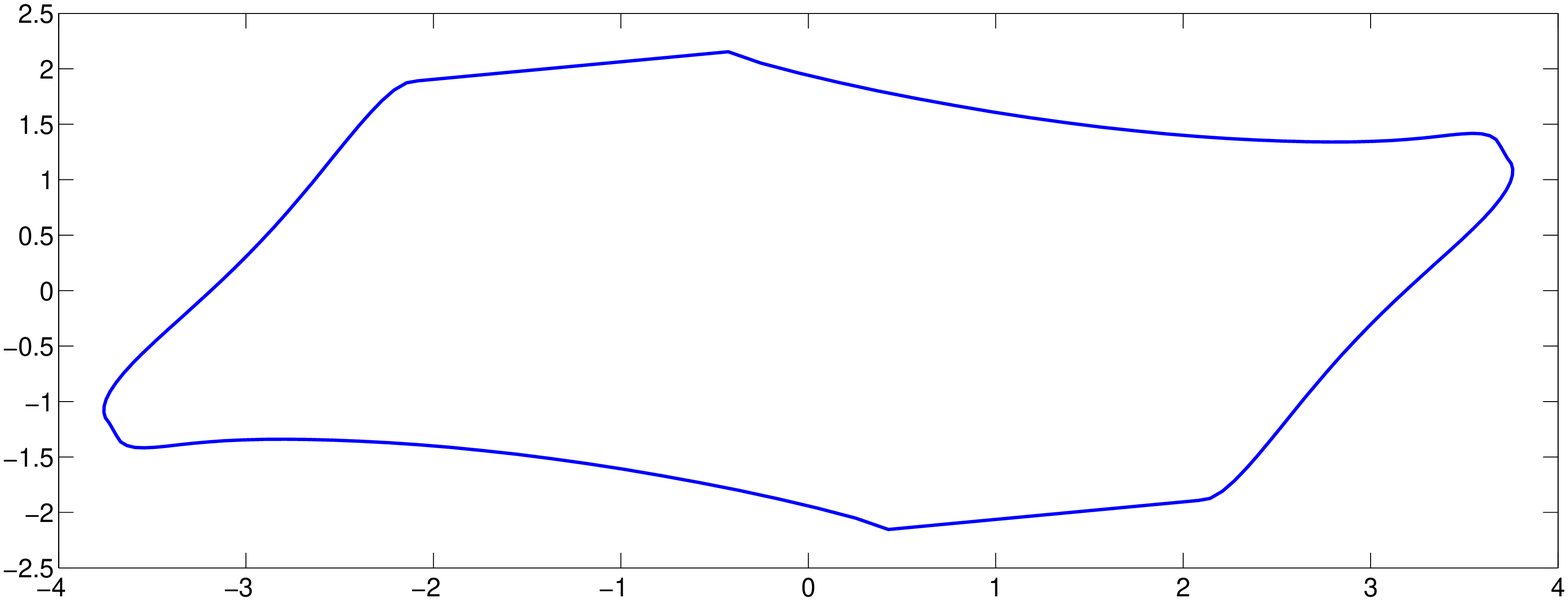}
  \includegraphics[height=45mm,width=60mm]{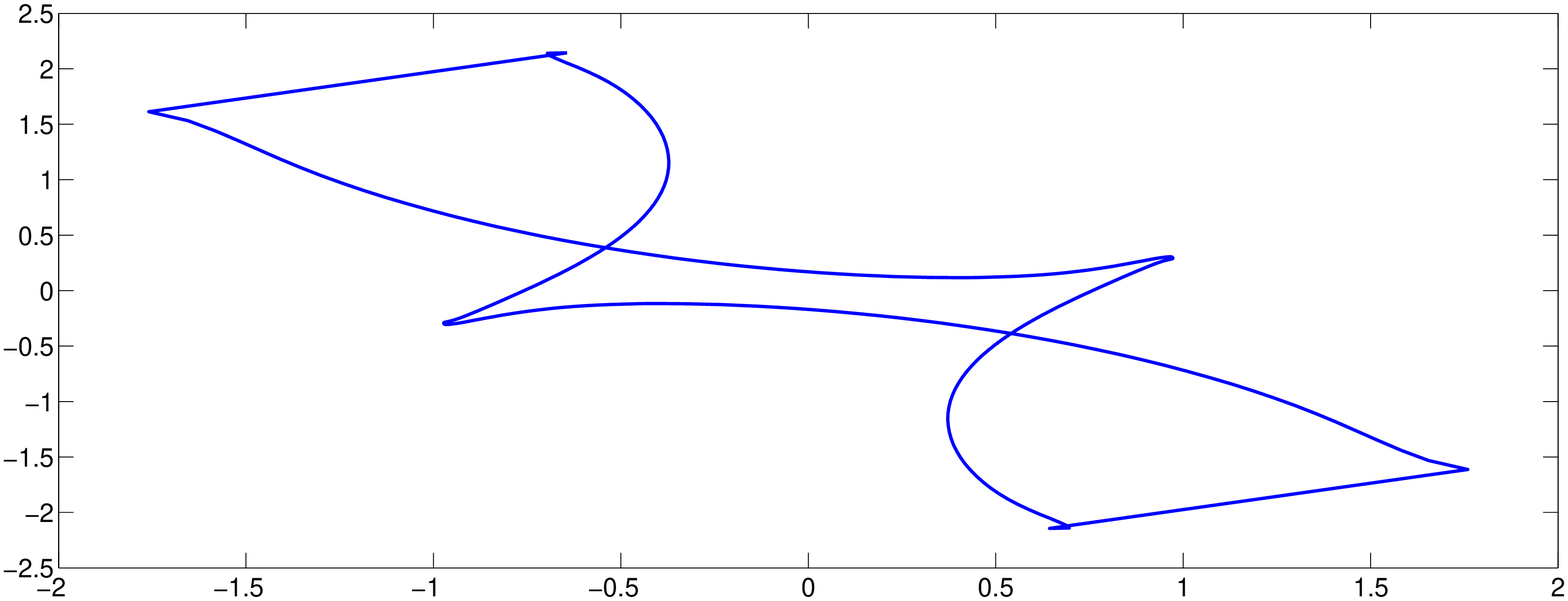}
\includegraphics[height=45mm,width=60mm]{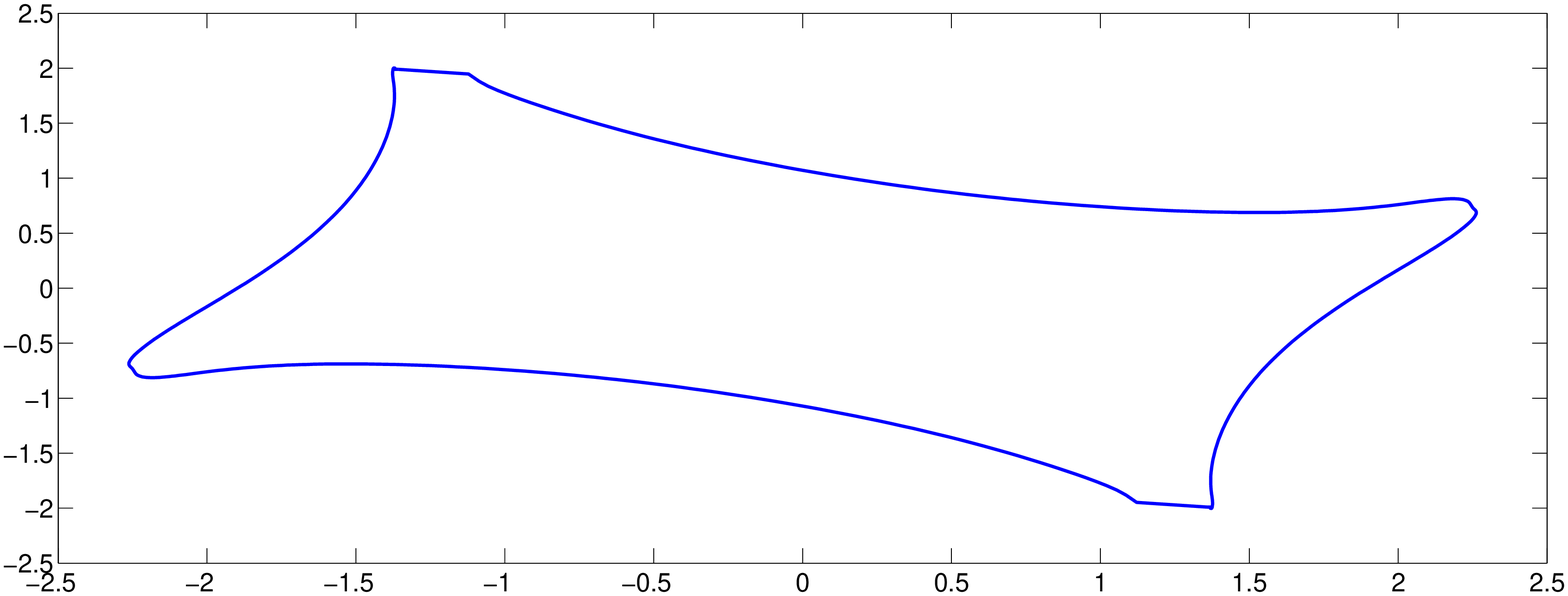}
  \includegraphics[height=45mm,width=60mm]{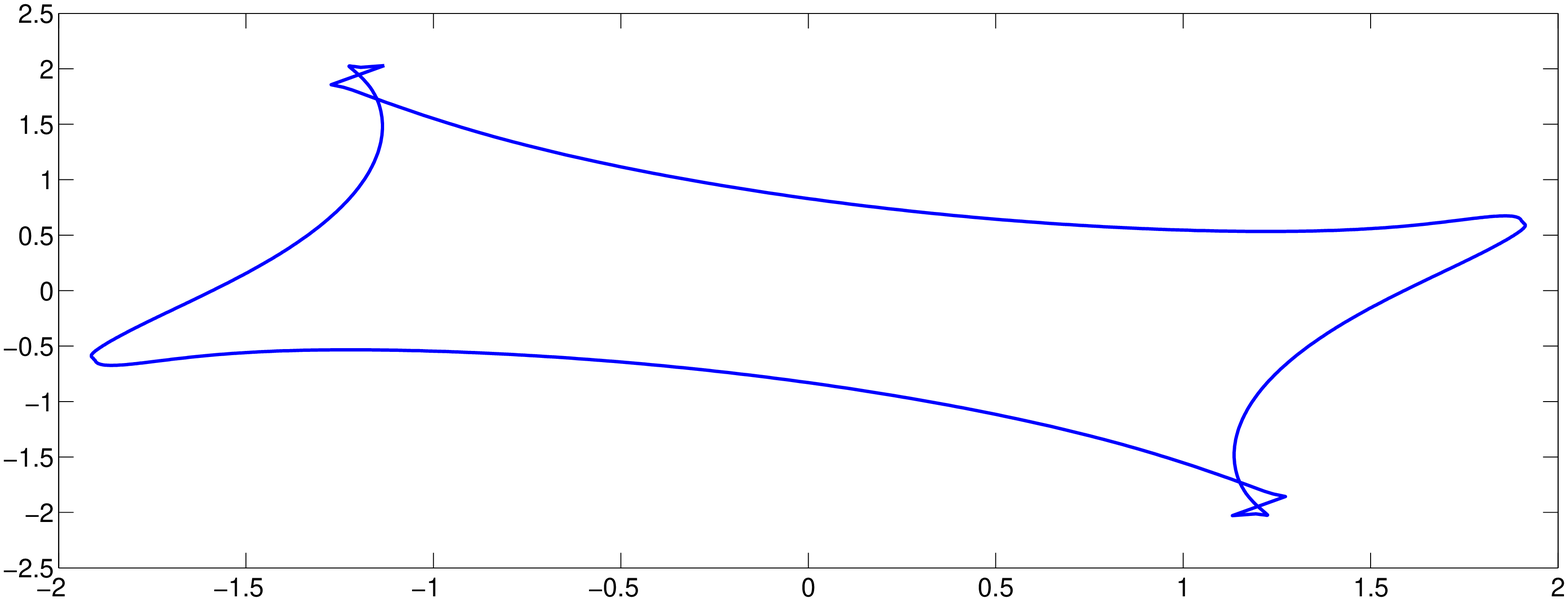}
\caption{\sf  Approximate solution $\omega_n(t)$ of the
Sherman--Lauricella equation \eqref{eq3} on the rhombus $\Gamma$,
$\alpha=\pi/5$ with $f:=f_1$ defined by \eqref{RHS} and $d=0$. From
the left to the right: $n=128,256,512,1024$} \label{fig:om_rh_pi5}
\end{figure}

Thus if the right hand sides $f\in W_2^1(\Gamma)$, the corrected
Sherman-Lauricella equation can be used in order to find an exact or
an approximate solution of the equation \eqref{eq3}. In the present
paper, we employ spline based Galerkin methods to the equation
\eqref{eq12} and study their stability and convergence. Let us
describe these methods in more detail. First of all, we normalize
all the basis spline functions used. If $n$ is fixed, then for any
$j\in \sZ$ the norm $||\widehat{\phi}_{nj}||$ of any basis element
$\widehat{\phi}_{nj}$ is
 $$
||\widehat{\phi}_{nj}||^2=\frac{1}{n}\int_0^d \widehat{\phi}^2(s)\,
ds.
 $$
Therefore, if $\nu_d$ refers to the number
\begin{equation}\label{cd}
\nu_d:= \left ( \int_0^d \widehat{\phi}^2(s)\, ds. \right )^{-1/2},
\end{equation}
then
\begin{equation}\label{eq12.5}
\phi_{nj}:= \nu_d \sqrt{n}\, \widehat{\phi}_{nj}, \quad j\in \sZ
\end{equation}
are unit norm vectors. An approximate solution of the equation
\eqref{eq12} is sought in the form
 \begin{equation}\label{eq13}
   \omega_n (t) =\sum_{\phi_{nk}\in S_n^d(\Gamma)} a_k \phi_{nk}(t),
\end{equation}
the coefficients $a_k$ of which are  obtained from the following
system of algebraic equations
 \begin{equation}\label{eq14}
    ( B_\Gamma \omega_n, \phi_{nj})=(f, \phi_{nj}), \quad \phi_{nj}\in
    S_n^d(\Gamma).
\end{equation}
An important problem now is to study the solvability of the
equations \eqref{eq14} and convergence of the approximate solutions
to an exact solution of the original Sherman--Lauricella equation
\eqref{eq3}. In Section \ref{s3}, this problem is discussed in a
more detail but, at the moment, we would like to illustrate the
method under consideration by a few numerical examples. Thus we
present Galerkin solutions of the equation \eqref{eq3} with the
right-hand side $f=f_1$,
 \begin{equation}\label{RHS}
 f_1(z) = f(x,y) = 4x^3-12xy^2+i(4y^3-12x^2y); \quad z=x+iy \in
 \Gamma,
 \end{equation}
on the unite square and rhombuses, and trace the evolution of the
solution when the initial contour is transformed from the unit
square into rhombuses with various opening angle $\alpha$. Some of
these contours have been used in \cite{DH:2011} in order to
illustrate the behaviour of the Nystr\"om method. Note that in the
corresponding examples from \cite{DH:2011}, approximate solutions of
the equation \eqref{eq3} with the right-hand side
 $$
f_2(z) = |z|
 $$
have been determined. We apply the spline Galerkin method to the
equations with such right-hand side, too. The results obtained have
a very good correlation with \cite{DH:2011} and the error evaluation
for both cases are reported in  Table \ref{table1}, where
$E^{f_i}_{n,\alpha}$ denotes the relative error $\| \omega _{2n}-
\omega _n\| _2/\| \omega _{2n}\| _2$ computed for the righthand side
$f_i$ and equation \eqref{eq3} is considered on the rhombus with the
opening angle $\alpha$.
 \begin{table}[!ht]
 \caption{Relative error of the spline Galerkin methods}
\begin{tabular}{|c|c|c|c|c||c|c|c|}
\hline n& $E^{f_1}_{n,{\pi}/{2}}$ & $E^{f_1}_{n,{\pi}/{3}}$ &
$E^{f_1}_{n,{\pi}/{4}}$ &
$E^{f_1}_{n,{\pi}/{5}}$&$E^{f_2}_{n,{\pi}/{2}}$ & $E^{f_2}_{n,{\pi}/{3}}$ & $E^{f_2}_{n,{\pi}/{6}}$ \\
\hline 128 & 0.0373 & 0.6194 & 1.3577  & 2.1716 & 0.0121 & 0.0217 & 0.0205\\
\hline 256 & 0.0198 & 0.0268 & 0.2046 & 0.6169 & 0.0067 & 0.0112 & 0.0245\\
\hline 512 & 0.0096 & 0.0059 & 0.0616 & 0.1888 & 0.0045 & 0.0102 & 0.0193 \\
\hline
\end{tabular}
  \label{table1}
\end{table}
In addition, Figures \ref{fig:square}--\ref{fig:om_rh_pi5} show the
convergence of the approximate solutions of the equation \eqref{eq3}
with the right-hand side \eqref{RHS} obtained by the Galerkin method
based on the splines of degree $d=0$ and the transformation of these
approximate solutions when $n$ increases.

Let us mention a few technical details related to the examples
below. Thus the rhombus with an opening angle $\alpha$ is
parameterized as follows,
\begin{equation}
\gamma (s) = \begin{cases}
\D 4s-\cos \left (\frac{\alpha}{2} \right )e^{i\alpha/2} & \text {
if }\, 0 \leq s
<1/4,\\[1ex]
\D (4s-1)e^{i\alpha}-i\sin \left (\frac{\alpha}{2} \right ) e^{i\alpha/2} & \text { if }\, 1/4 \leq s < 1/2,\\[1ex]
\D -(4s-2)+\cos \left (\frac{\alpha}{2} \right )e^{i\alpha/2} & \text { if }\, 1/2 \leq s < 3/4,\\[1ex]
\D -(4s-3)e^{i\alpha} + i\sin \left (\frac{\alpha}{2} \right
)e^{i\alpha/2}&\text { if }\, 3/4 \leq s \leq 1.
\end{cases}
\end{equation}
Moreover, we have to compute the scalar products $(B_\Gamma \omega
_n, \phi _{nj}) $.
 Recall that $\supp \phi _{nj} \subset [j/n, \;
(j+d+1)/n]$ and use the Gauss-Legendre quadrature rule with
quadrature points which coincide with the zeros of the Legendre
polynomial $P_{24}(x)$ on the canonical interval $[-1,1]$, scaled
and shifted to the interval $[j/n, \; (j+d+1)/n]$. More
specifically, the corresponding formula is
  \begin{equation}\label{quad1}
 (B_\Gamma \omega _n, \phi _{nj})=\int_{j/n}^{(j+d+1)/n} B_\Gamma
\omega _n(\phi (s)) \overline{\phi _{nj}(\phi (s))}ds\approx
\sum_{k=1}^{24}w_kB_\Gamma \omega _n(\phi (s_k)) \phi _{nj}(\phi
(s_k)),
 \end{equation}
where $w_k,\, s_k$ are the Gauss-Legendre weights and the
Gauss-Legendre points on the interval $[j/n, \; (j+d+1)/n]$. In
order to find the values of the corresponding line integrals at the
Gauss-Legendre points, the composite Gauss-Legendre quadrature is
used \cite[Section 3]{DH:2011},  namely,
 \begin{equation}\label{quad2}
 \begin{aligned}
\ig k(t,\tau ) x(\tau ) d\tau & = \int_0^1 k(\gamma (\sigma ),
\gamma (s))x(\gamma (s))\gamma ' (s) ds \\
& \approx \sum_{l=0}^{m-1} \sum \limits _{p=0}^{r-1}w_p
k(\gamma (\sigma ),  \gamma (s_{lp}))x(\gamma (s_{lp}))\tau '_{lp}/m,\\
\end{aligned}
 \end{equation}
where $\tau '_{lp} = \gamma '(s_{lp})$ with $m=40$ and $r=24$.

Table \ref{table1} and Figures
\ref{fig:square}-\ref{fig:om_rh_unit_pi3} show a good convergence of
approximate solutions if the corner point of the contour has an
opening angle close or equal to $\pi/2$. On the other hand, the
presence of opening angles of a small magnitude can cause problems
and lead to a convergence slowdown (see Figures
\ref{fig:om_rh_pi4}-\ref{fig:om_rh_pi5}). Note that although the
focus of this work is on the stability, the error estimates
presented in Table \ref{table1}  are comparable with estimates of
recent work \cite{JWB:2014} for fast Fourier--Galerkin method for an
integral equation used to solve boundary value problem \eqref{eq2}
in smooth domains. Moreover, further improvement of the convergence
rate is possible if for the approximations of singular integrals and
inner products arising in the Galerkin method one employs graded
meshes of various kind \cite{CTT:2004, He:2011}.

\section{Galerkin method. Local operators and stability\label{s3}}

Our next task is to find conditions of applicability of the spline
Galerkin methods to the equation \eqref{eq12}. It is worth
mentioning that for smooth contours $\Gamma$, the methods considered
here are always applicable and provide satisfactory results. For
details the reader can consult \cite{DV:2007}, where similar methods
for the Muskhelishvili equation on smooth contours are considered.
On the other hand, the presence of corners changes the situation
drastically, and the applicability of the approximation method is
not always guaranteed.

Let $P_n$ be the orthogonal projection from $L^2(\Gamma)$ on the
subspace $S_n^d(\Gamma)$. Then the systems \eqref{eq14} is
equivalent to the following operator equations
 \begin{equation}\label{eq15}
    P_n B_\Gamma P_n \omega_n = P_n f, \quad n\in \sN.
\end{equation}

 \begin{definition}
We say that the sequence $( P_n B_\Gamma P_n)$ is stable if there is
an $m\in \sR$ and an $n_0\in\sN$ such that for all $n\geq n_0$ the
operators $P_n B_\Gamma P_n :S_n^d(\Gamma) \to S_n^d(\Gamma)$ are
invertible and
 $$
|| (P_n B_\Gamma P_n)^{-1}P_n|| \leq m
 $$
for all $n\geq n_0$.
 \end{definition}

Recall that if the stability of the corresponding sequence $(P_n
B_\Gamma P_n)$ is established, then the convergence of the Galerkin
method and error estimates can be obtained from well known results,
cf. \cite[Section 1.6, inequality (1.30)]{DS:2008}. Therefore, in
this work we mainly deal with the stability and our approach is
based on $C^*$-algebra methods often used in operator theory. Let
$\cL_{add}(L^2(\Gamma))$ refer to the real $C^*$-algebra of all
additive continuous operators on the space $L^2(\Gamma)$. One can
show \cite{DS:2008} that any operator $A\in \cL_{add}(L^2(\Gamma))$
admits the unique representation $A=A_1 + A_2 M$ where $A_1, A_2$
are linear operators and $M$ is the operator of complex conjugation.
This representation allows one to introduce the operation of
involution on $\cL_{add}(L^2(\Gamma))$ as follows
\begin{equation}\label{eq17}
A^*:=A_1^* + M A_2^*,
\end{equation}
with $A_1^*,A_2^*$ being usual adjoint operators to the linear
operators $A_1,A_2$,  cf. \cite[Theorem 1.3.8 and Example
1.3.9]{DS:2008}. By $\cA^\Gamma$ we denote the set of all bounded
sequences $(A_n)$ of bounded additive operators $A_n:\im P_n\to \im
P_n$ such that there is an operator $A\in \cL_{add}(L^2(\Gamma))$
with the property
\begin{equation*}
   s-\lim A_n P_n= A, \quad s-\lim\, (A_n P_n)^*P_n= A^*,
\end{equation*}
where $s-\lim A_n$ denotes the strong limit of the operator sequence
$(A_n)$.

Provided with natural operations of addition, multiplication,
multiplication by scalars $\lambda\in \sC$, with an involution
introduced according to \eqref{eq17}, and with the norm
 $$
||(A_n)||:=\sup_{n\in\sN}||A_n||,
 $$
the set $\cA^\Gamma$ becomes a real $C^*$-algebra. Consider also the
subset $\cJ^\Gamma\subset \cA^\Gamma$ consisting of all sequences
$(J_n)$ of operators $J_n:\im P_n\to \im P_n$ which can be
represented in the form
 $$
J_n=P_n T P_n + C_n, \quad n\in \sN,
 $$
where the operator $T$ belongs to the ideal
$\cK_{add}(L^2(\Gamma))\subset \cL_{add}(L^2(\Gamma))$ of all
compact operators and the sequence $(C_n)$ tends to zero uniformly,
i.e.
 $$
\lim_{n\to \infty}||C_n||=0.
 $$
The stability of sequences from the algebra $\cA^\Gamma$ can be
characterized as follows.

 \begin{theorem}[cf. {\cite[Proposition 1.6.3]{DS:2008}}]\label{t2}
A sequence $(A_n)\in \cA^\Gamma$ such that $A:=s-\lim A_nP_n$ is
stable if and only if the operator $A$ is invertible in
$\cL_{add}(L^2(\Gamma))$ and the coset $(A_n)+\cJ^\Gamma$ is
invertible in the quotient algebra $\cA^\Gamma/\cJ^\Gamma$.
 \end{theorem}

Consider now the sequence $(P_n B_\Gamma P_n)$ of the Galerkin
operators defined by the projection operators $P_n$. Recall that on
the space $L^2(\Gamma)$ the sequence of the orthogonal projections
$(P_n)$ strongly converges to the identity operator $I$ and
$P_n^*=P_n, n\in \sN$. It implies that for any operator $A\in
\cL_{add}(L^2(\Gamma))$ the following relations
 $$
s-\lim P_n A P_n = A, \quad s-\lim (P_n A P_n)^*P_n = A^*
 $$
hold \cite{PS:1991}. Therefore, combining Theorem \ref{t1}\/ and
Theorem \ref{t2} one obtains the following result.

 \begin{corollary}\label{c1}
Let $\Gamma$ be a simple closed piecewise smooth curve. The spline
Ga\-ler\-kin method \eqref{eq15} is stable if and only  if the coset
$(P_n B_\Gamma P_n) +\cJ^\Gamma$ is invertible in the quotient
algebra $\cA^\Gamma/\cJ^\Gamma$.
 \end{corollary}

Thus in order to establish the stability of the Galerkin method, one
has to study the invertibility of the coset $(P_n B_\Gamma P_n)
+\cJ^\Gamma$ in the algebra $\cA^\Gamma/\cJ^\Gamma$. This problem
can be tackled more efficiently, if we restrict ourselves to a
smaller algebra containing the coset $(P_n B_\Gamma P_n) +\cJ$. More
precisely, let $M$ refer to the operator of the complex conjugation,
 $$
 M\phi(t):=\overline{\phi(t)}, \quad t\in \Gamma,
 $$
and let $S_\Gamma$ be the Cauchy singular integral operator,
 $$
S_\Gamma\phi(t):=\frac{1}{\pi i}\int_\Gamma \frac{\phi(\zeta)}{\zeta
- t}\,d\zeta.
 $$
Consider the smallest closed real $C^*$-subalgebra  $\cB^\Gamma$ of
the algebra $\cA^\Gamma$ which contains all operator sequences of
the form $(P_n M P_n), (P_n S_\Gamma P_n)$ and also the sequences
$(P_n f P_n), f\in \sC_\sR(\Gamma)$ and $(G_n)$, where
$\lim_{n\to\infty} ||G_n||=0$ and $\sC_\sR(\Gamma)$ is the set of
all  continuous real-valued functions on the contour $\Gamma$.

 \begin{remark}\label{r1}
It follows from  \cite{DS:2002, Du:1979, Mikhlin1986, PS:1991} that
$\cJ^\Gamma \subset \cB^\Gamma$ and that the sequence $(P_n B_\Gamma
P_n)$ belongs to $\cB^\Gamma$. Therefore, $\cB^\Gamma/\cJ^\Gamma$ is
a real $C^*$-subalgebra of $\cA^\Gamma/\cJ^\Gamma$, and  by
\cite[Corollary 1.4.10]{DS:2008} the coset $(P_n B_\Gamma
P_n)+\cJ^\Gamma$ is invertible in $\cA^\Gamma/\cJ^\Gamma$ if and
only if it is invertible in $\cB^\Gamma/\cJ^\Gamma$. Therefore, one
can now study the invertibility of the coset $(P_n B_\Gamma P_n)
+\cJ^\Gamma$ in the smaller algebra $\cB^\Gamma/\cJ^\Gamma$. To this
end we will employ a localizing principle.
 \end{remark}

Thus with each point $\tau\in\Gamma$ we associate a model contour
$\Gamma_\tau$ as follows. Let $\theta_\tau$ be the angle between the
right and the left semi-tangents to $\Gamma$ at the point $\tau$,
and let $\beta_\tau$ refer to the angle between the right
semi-tangent to $\Gamma$ and the real line $\sR$. Consider now the
curve
  $$
\Gamma_\tau:=e^{i(\beta_\tau+\theta_\tau)} \sR_+^- \cup
e^{i\beta_\tau} \sR_+^+
  $$
where $\sR_-^+$ and $\sR_+^+$ denote the positive semi-axis $\sR^+$
correspondingly directed to and out of the origin. Further, on each
such contour $\Gamma_\tau$, $\tau\in \Gamma$ we consider the
corresponding Sherman-Lauricella operator
 \begin{equation}\label{eq16}
A_\tau=I + L_\tau-K_\tau M,
\end{equation}
where
 \begin{equation*}
    L_\tau\ot :=\frac{1}{2\pi i} \int_{\Gamma_{\tau}} \o \, d\ln \left ( \frac{\zeta-t}{\ao
    -\tov}\right ), \quad
    K_\tau \ot:= \frac{1}{2\pi i} \int_{\Gamma_{\tau}} \o \, d \left ( \frac{\ao
    -\tov}{\zeta-t}\right ).
 \end{equation*}

Analogously to the algebra $\cB^\Gamma$ and to the ideal
$\cJ^\Gamma$ one can introduce algebras $\cB^{\Gamma_\tau}$ and
ideals $\cJ^{\Gamma_\tau}\subset \cB^{\Gamma_\tau}$, $\tau\in
\Gamma$, which allow to establish conditions of the applicability of
the corresponding Galerkin method for the operators \eqref{eq16}.
For this we also need appropriate spline spaces on both the contour
$\Gamma_\tau$ and the positive semi-axis $\sR^+:=\sR^+_+$. These
spline spaces can be constructed by using the functions
\eqref{eq12.5} again. More precisely, consider the functions
  \begin{equation}\label{eq17.5}
\widetilde{\phi}_{nj}(t):=\left \{
 \begin{array}{l}
 \left \{
 \begin{array}{cc}
   \phi_{nj}(s) & \D \text{if} \quad t=e^{i\beta_\tau}s \\
   0 & \text{otherwise} \\
 \end{array}
   \right .
   \quad j\geq 0\\[3ex]
    \left \{
\begin{array}{cc}
   \phi_{n,j-d}(s) & \D \text{if} \quad t=e^{i(\beta_\tau+\theta_\tau)}s \\
   0 & \text{otherwise} \\
 \end{array}
 \right.
      \quad j< 0
 \end{array}.
   \right.
\end{equation}
Let $S_n^d(\Gamma_\tau)$  and $S_n^d(\sR^+)$ be, respectively, the
smallest closed subspaces of $L_2(\Gamma_{\tau})$ and $L^2(\sR^+)$
which contains all functions \eqref{eq17.5} and  all functions
$\phi_{nj}$, $j\geq0$ of \eqref{eq17.5} for $\beta_\tau =0$.
Moreover, let $P_n^\tau$, $n\in\sN$ and $P_n^+$ denote the
orthogonal projection onto subspaces $S_n^d(\Gamma_\tau)$  and
$S_n^d(\sR^+)$, respectively. In order to study the stability of the
sequence $(P_n^\tau A_\tau P_n^\tau )$, one can apply Theorem
\ref{t2} and Remark \ref{r1} to obtain the following result.

\begin{corollary}\label{c2}
The sequence $(P_n^\tau A_\tau P_n^\tau )\in \cB^{\Gamma_\tau}$  is
stable if and only if the operator $A_\tau$ is invertible in
$\cB^\tau$ and the coset $(P_n^\tau A_\tau P_n^\tau
)+\cJ^{\Gamma_\tau}$ is invertible in the quotient algebra
$\cB^{\Gamma_\tau}/\cJ^{\Gamma_\tau}$.
 \end{corollary}

Further, let $L^2_2(\sR^+)$ be the space  of all pairs
$(\varphi_1,\varphi_2)^T, \, \varphi_1,\varphi_2\in L^2(\sR^+) $
provided with the norm
 $$
||(\varphi_1,\varphi_2)^T||:= \left ( \right
||\varphi_1||^2+||\varphi_2||^2)^{1/2},
 $$
and let $\eta : L^2(\Gamma_\tau)\to L^2_2(\sR^+)$ be the mapping
defined by
 $$
\eta (\vp)=(\vp(se^{i(\beta_\tau +
\omega_\tau)}),\vp(se^{i\beta}))^T, \quad s\in \mathbb{R}^+,
 $$
where $a^T$ denotes the transposition of the vector $a$. It is clear
that $\eta$ is a linear isometry from $L^2(\Gamma_\tau)$ onto
$L^2_2(\sR^+)$. Moreover, the mapping
$\Psi:\cL_{add}(L^2(\Gamma_\tau))\to \cL_{add} (L^2_2(\sR^+))$
defined by
 \begin{equation}\label{EqPsi}
\Psi(A)=\eta A\eta^{-1},
 \end{equation}
is an isometric algebra isomorphism. In particular,  straightforward
calculations show that

\begin{align}
\Psi(P_n^\tau)&=\diag (P_n^+, P_n^+) , \label{eq18}\\[1.5ex]
  \Psi(M)&=\diag(\widetilde{M},\widetilde{M}), \label{eq18.5}\\[1.5ex]
   \Psi(L_\tau) &  =\left (
 \begin{array}{cc}
  0  & \mathcal{N}_{\theta_\tau} \\[1ex]
 \mathcal{N}_{\theta_\tau}    & 0 \\
    \end{array}
    \right )    ,           \label{eq19}\\[1.5ex]
    \Psi(K_\tau)& = \left (
 \begin{array}{cc}
0    & e^{i2\beta_\tau}\mathcal{M}_{2\pi-\theta_\tau} \\[1ex]
  -e^{i2(\beta_\tau+\theta_\tau)}\mathcal{M}_{\theta_{j}} & 0  \\
    \end{array}
    \right ) ,   \label{eq20}
\end{align}
 where
\begin{align*}
 &\mathcal{N}_{\theta_\tau} \varphi(\sigma)=\frac{1}{2}\frac{1}{2\pi i}\int_0^{\infty}
 \left (  \frac{1}{1-(\sigma/s) e^{i\theta_\tau}} -
\frac{1}{1-(\sigma/s) e^{i(2\pi-\theta_\tau)}}  \right )\varphi(s)\frac{ds}{s} ,\\[1.5ex]
   & \cM_{\theta_\tau} \varphi(\sigma):=\frac{1}{\pi}\int_{0}^{\infty}\left (
\frac{\sigma}{s}\right ) \frac{\sin
\theta_\tau}{(1-(\sigma/s)e^{i\theta_\tau})^2} \varphi(s)\frac{ds}{s} , 
\end{align*}
and the symbol $\widetilde{M}$ in the right-hand  side of
\eqref{eq18.5} refers to the operator of the complex conjugation on
the space $L^2(\sR^+)$. Moreover, one can observe that the operators
$\mathcal{N}_{\theta_\tau}$ and $\mathcal{M}_{\theta_\tau}$ have a
special form -- viz.
\begin{equation}\label{eq19.5}
K \varphi (\sigma):=\int_0^\infty
\mathsf{k}_{\theta_\tau}\left(\frac{\sigma}{s} \right ) \varphi(s)
\frac{ds}{s}
\end{equation}
and
  \begin{align}
   \mathsf{k}_{\theta_\tau}=\mathsf{k}_{\theta_\tau}(u)&:=\mathsf{n}_{\theta_\tau}(u)=\frac{1}{2\pi}
   \frac{u\sin\theta_\tau}{|1-ue^{i\theta_\tau}|^2}
   ,\quad \text{if} \quad K=\mathcal{N}_{\theta_\tau},\label{eq21}\\[1ex]
 \mathsf{k}_{\theta_\tau}=\mathsf{k}_{\theta_\tau}(u)&:=\mathsf{m}_{\theta_\tau}(u)=\frac{1}{\pi}
   \frac{u\sin\theta_\tau}{(1-ue^{i\theta_\tau})^2}
   ,\quad\text{if}\quad K=\mathcal{M}_{\theta_\tau}.\label{eq22}
\end{align}
On the space $l^2$ of the sequences $(\xi_k)$ of complex numbers
$\xi_k, k=0,1,\ldots$,
 $$
l^2:=\{ (\xi_k)_{k=0}^\infty: \sum_{k=0}^\infty |\xi_k|^2 <\infty
\},
 $$
the function $\textsf{k}_{\theta_\tau}$ defines a bounded linear
operator $A(\mathsf{k}_{\theta_\tau})$ with the matrix
representation
 $$
A(\mathsf{k}_{\theta_\tau})= \left ( \nu_d^2 \int_0^{d+1}
\widehat{\phi}(t) \int_0^{d+1} \mathsf{k}_{\theta_\tau} \left
(\frac{u+l}{t+q}\right ) \widehat{\phi}(u) \frac{du}{u+q} \,dt\right
)_{q,l=0}^\infty\, ,
 $$
where $\nu_d$ is the constant \eqref{cd}.

 \begin{theorem}\label{t3}
Let $\mathsf{n_{\theta_\tau}}$ and $\mathsf{m_{\theta_\tau}}$ be the
functions defined  by \eqref{eq21} and \eqref{eq22}, respectively.
The spline Galerkin method \eqref{eq15} is stable if and only if the
operators $R^\tau:l^2\times l^2 \to l^2\times l^2$,%
\begin{equation}\label{eq23}
 \begin{aligned}
  & R^\tau:= \\[1ex]
   &  \left(%
\begin{array}{c@{\hspace{-1mm}}c}
 I  & A(\mathsf{n}_{\theta_\tau}) \\[1ex]
   A(\mathsf{n}_{\theta_\tau})  & I \\
   \end{array}%
\right) \!\! +\!\!
 \left(%
\begin{array}{c@{\hspace{-1mm}}c}
 0  &  e^{i\beta_\tau} A(\mathsf{m}_{2\pi-\theta_\tau})\\[1ex]
  -e^{-i(\beta_\tau+\theta_\tau)} A(\mathsf{m}_{\theta_\tau})  & 0 \\
   \end{array}%
\right)
 \left(%
\begin{array}{cc}
 \overline{M}  & 0 \\[1ex]
  0  &  \overline{M}\\
   \end{array}%
\right)
\end{aligned}
\end{equation}
are invertible for all $\tau\in \cM_\Gamma$.
 \end{theorem}

 \textbf{Proof.}
By Corollary \ref{c1} the sequence $(P_n B_\Gamma P_n)$ is stable if
and only if the coset $(P_n B_\Gamma P_n) +\cJ^\Gamma$ is
invertible. Moreover, since  $T_\Gamma$ of \eqref{eq10} is a compact
operator, the sequences $(P_n A_\Gamma  P_n)$ and $(P_n B_\Gamma
P_n)$ belong to the same coset $(P_n A_\Gamma  P_n)+\cJ^\Gamma$ of
the quotient algebra  $\cB^\Gamma/\cJ^\Gamma$. However, by a version
of the Allan's Local Principle \cite{Allan1968} for real
$C^*$-algebras \cite[Theorem 1.9.5]{DS:2008}, the coset $(P_n
A_\Gamma P_n)+\cJ^\Gamma$ is invertible if and only if for every
$\tau\in\Gamma$ the coset $(P_n^\tau A_{\Gamma_\tau}
P_n^\tau)+\cJ^{\Gamma_\tau}$ is invertible in the corresponding
algebra $\cB^{\Gamma_\tau}/ \cJ^{\Gamma_\tau}$. Therefore, the
stability of our operator sequence will be established if we manage
to show the invertibility of all cosets $(P_n^\tau A_{\Gamma_\tau}
P_n^\tau)+\cJ^{\Gamma_\tau}$, $\tau\in \Gamma$. Let us start with
the case where $\tau$ is not a corner point of $\Gamma$. If
$\tau\notin \cM_\Gamma$, then $\theta_\tau=\pi$, and straightforward
calculations show that $L_\tau$ and $K_\tau$ are the zero operators.
Hence, $A_\tau$ is just the identity operator $I$ in the
corresponding space, so that $P_n^\tau A_\tau P_n^\tau =P_n^\tau$.
The sequence $(P_n^\tau)$ is obviously stable so that the
corresponding coset $(P_n^\tau) +\cJ^\tau$ is invertible.

Consider next the case where $\tau\in\cM_\Gamma$. Note that by
\cite[Theorem 2.2]{DH:2011} the operator $A_\tau$ is invertible on
the space $L^2(\Gamma)$. Therefore, by Corollary \ref{c2} the coset
$(P_n^\tau A_\tau P_n^\tau) +\cJ^\tau$ is invertible in
$\cB^{\Gamma_\tau}/ \cJ^{\Gamma_\tau}$ if and only if the sequence
$(P_n^\tau A_\tau P_n^\tau)$ is stable. However, the stability of
this sequence is equivalent to the stability of the sequence
$(\Psi(P_n^\tau A_\tau P_n^\tau))$, where mapping $\Psi$ is defined
by \eqref{EqPsi}. Consider also the operators $\Lambda_n:
S_n^d(\sR^+)\to l^2$ defined by
  $$
\Lambda_n \left (\sum_{j=0}^\infty \xi_j \phi_{nj} \right )= (\xi_0,
\xi_1, \ldots, \,).
  $$
By Lemma \ref{l1} these operators are bounded and continuously
invertible. Set $\Lambda_{-n}:=\Lambda_n^{-1}$ and note that the
sequence $(\Psi(P_n^\tau A_\tau P_n^\tau))$ is stable if and only if
so is the sequence $(R_n^\tau)$, where
 $$
R_n^\tau=\diag(\Lambda_{n},\Lambda_{n})\cdot \Psi(P_n^\tau A_\tau
P_n^\tau)\cdot \diag(\Lambda_{-n},\Lambda_{-n} ): l^2 \times l^2\to
l^2 \times l^2.
 $$
From the definition of the mappings $\Psi$ and $\Lambda_{\pm n}$ one
obtains that the operators $R_n^\tau$ have the form
 \begin{equation*}
    R_n^\tau = (A_{lp}^{(n,\tau)})_{l,p=1}^2
    +(B_{lp}^{(n,\tau)})_{l,p=1}^2 \diag
    (\overline{M},\overline{M}),
\end{equation*}
with the operators $A_{lp}^{(n,\tau)},B_{lp}^{(n,\tau)}:l^2\to l^2$
defined according to the relations \eqref{eq18}-\eqref{eq20}.
However, these operators do not depend on the parameter $n$ at all.
Really, consider the matrix representations of the operators
$A_{12}^{(n,\tau)}, A_{21}^{(n,\tau)},B_{12}^{(n,\tau)},
B_{21}^{(n,\tau)}$.  It follows from \eqref{eq19.5} that the entries
$a_{lq}$ of the corresponding matrices $(a_{lq})_{l,q=0}^\infty$ are
 \begin{align*}
 a_{pq}& =\int_{\sR^+} K \phi_{qn}(\sigma)\phi_{ln}(\sigma)\,d\sigma
 =\int_{\sR^+} \int_{\sR^+}\mathsf{k}_{\theta_\tau} \left(\frac{\sigma}{s}\right)
 \phi(ns-q) \frac{ds}{s}\phi(n\sigma-l)\,d\sigma \\[1.5ex]
 &=\frac{1}{n}\int_{\sR^+} \int_{\sR^+}\mathsf{k}_{\theta_\tau} \left(\frac{u+l}{t+q}\right)
 \phi(u) \frac{du}{u+q}\phi(t)\,dt \\[1.5ex]
&=\frac{1}{n}\int_{\sR^+} \int_{\sR^+}\mathsf{k}_{\theta_\tau}
\left(\frac{u+l}{t+q}\right)
( \nu_d \sqrt{n}\widehat{\phi}(u)) \frac{du}{u+q}(\nu_d\sqrt{n}\widehat{\phi}(t))\,dt \\[1.5ex]
  &=\nu_d^2 \int_0^{d+1}
\widehat{\phi}(t) \int_0^{d+1} \mathsf{k}_{\theta_\tau} \left
(\frac{u+l}{t+q}\right ) \widehat{\phi}(u) \frac{du}{u+q} \,dt,
\end{align*}
hence these operators are independent of $n$. Moreover,
$B_{11}^{(n,\tau)}, B_{22}^{(n,\tau)}=0$ and $A_{11}^{(n,\tau)}=
A_{22}^{(n,\tau)}=I$. Combining all the above representations, one
obtains that the operators $R_n^\tau$ do not depend on the parameter
$n$. Therefore, $(R_n^\tau)$ is a constant sequence and it is stable
if and only if any its member, say $R^\tau_1$, is invertible. It
remains to observe that $R^\tau=R^\tau_1$, which completes the
proof.
  \rbx

\section{Numerical approach to the invertibility of local
operators}.

As was already mentioned, there is no efficient analytic method to
verify the invertibility of the local operators $R^\tau$. On the
other hand, numerical approaches turn out to be surprisingly
fruitful. Recall that the operators $R^\tau,\; \tau \in \cM_\Gamma$
do not depend on the shape of the contour $\Gamma$ but only on the
relevant angles $\theta_\tau$ and $\beta_\tau$. Therefore, for
contours having only one corner point, Theorem \ref{t3} can be
reformulated as follows.

 \begin{corollary}\label{c3}
If $\tau$ is the only corner point of the contour $\Gamma$, then the
operator $R^\tau$ is invertible if and only if the Galerkin method
$(P_n B_\Gamma P_n)$ is stable.
 \end{corollary}

\begin{figure}[!ht]
\centering
\includegraphics[height=45mm,width=60mm]{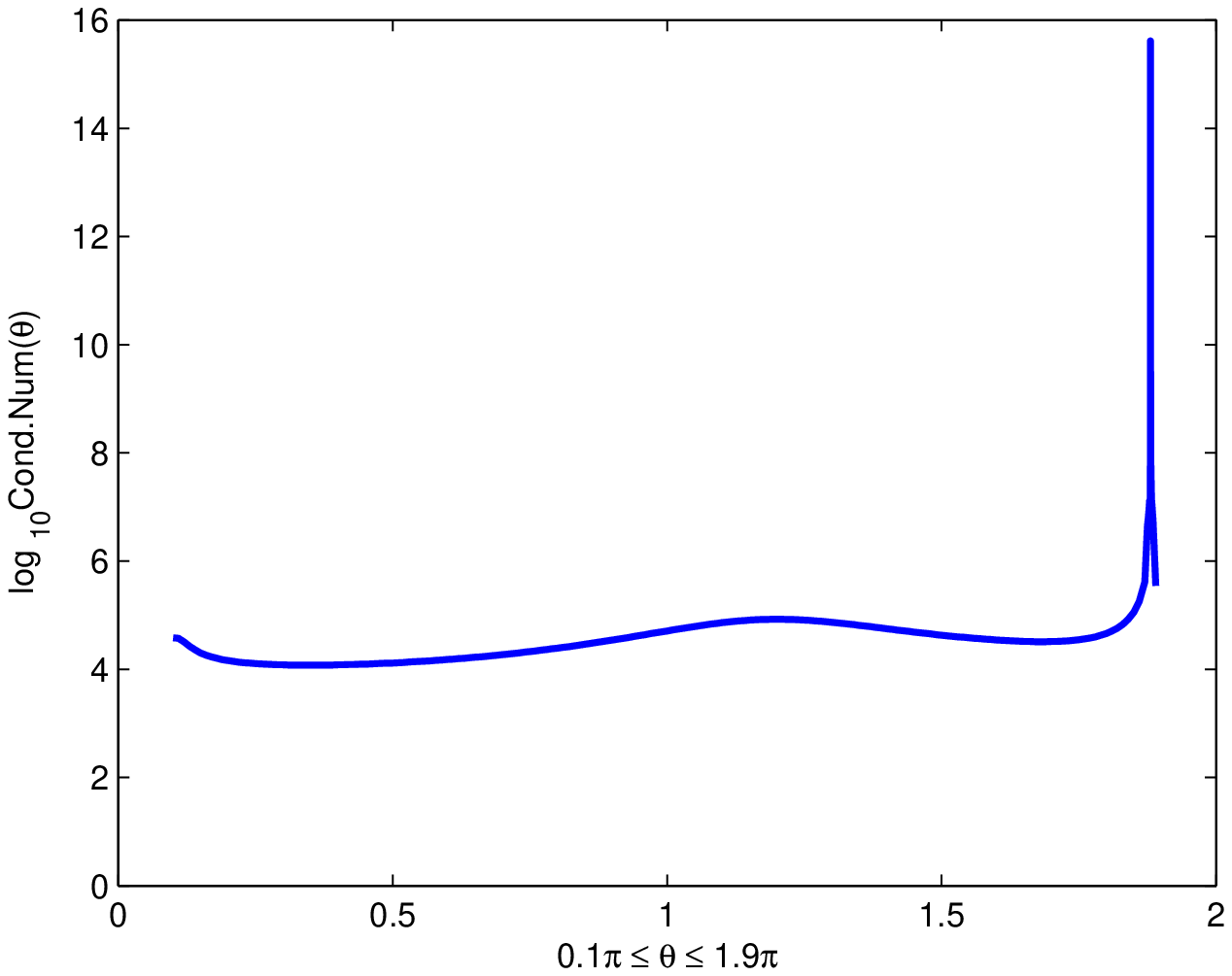}
\includegraphics[height=45mm,width=60mm]{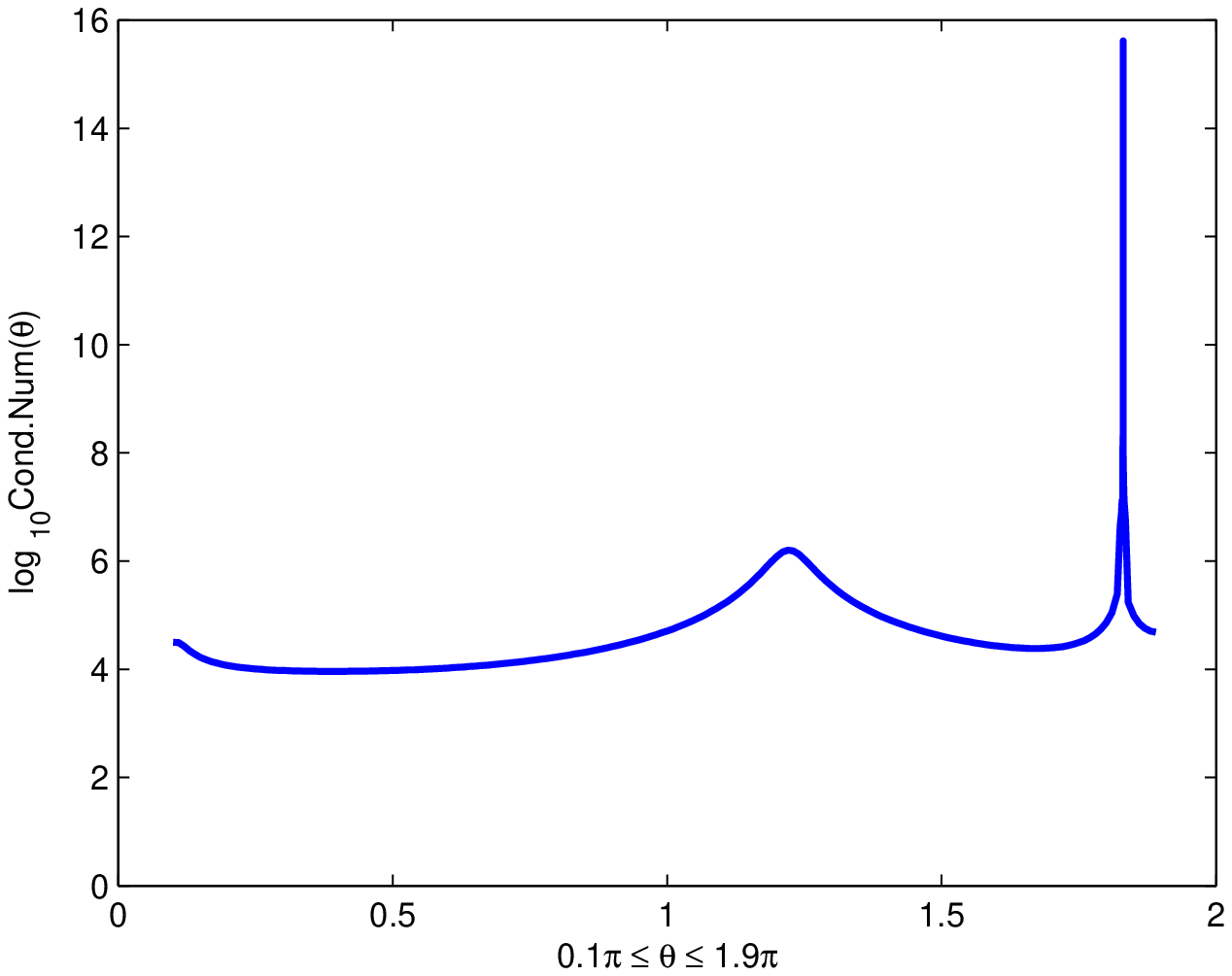}
\includegraphics[height=45mm,width=60mm]{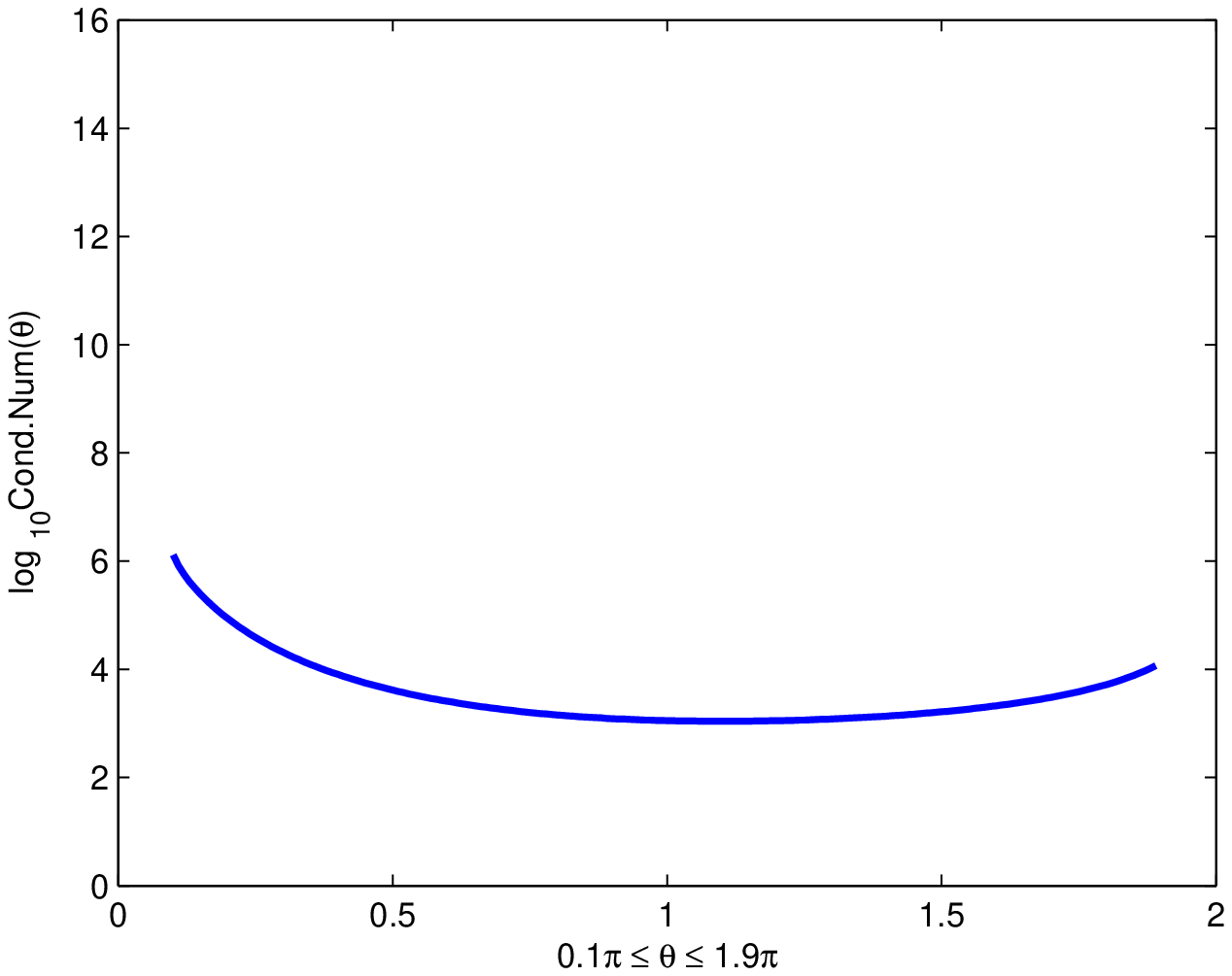}
\includegraphics[height=45mm,width=60mm]{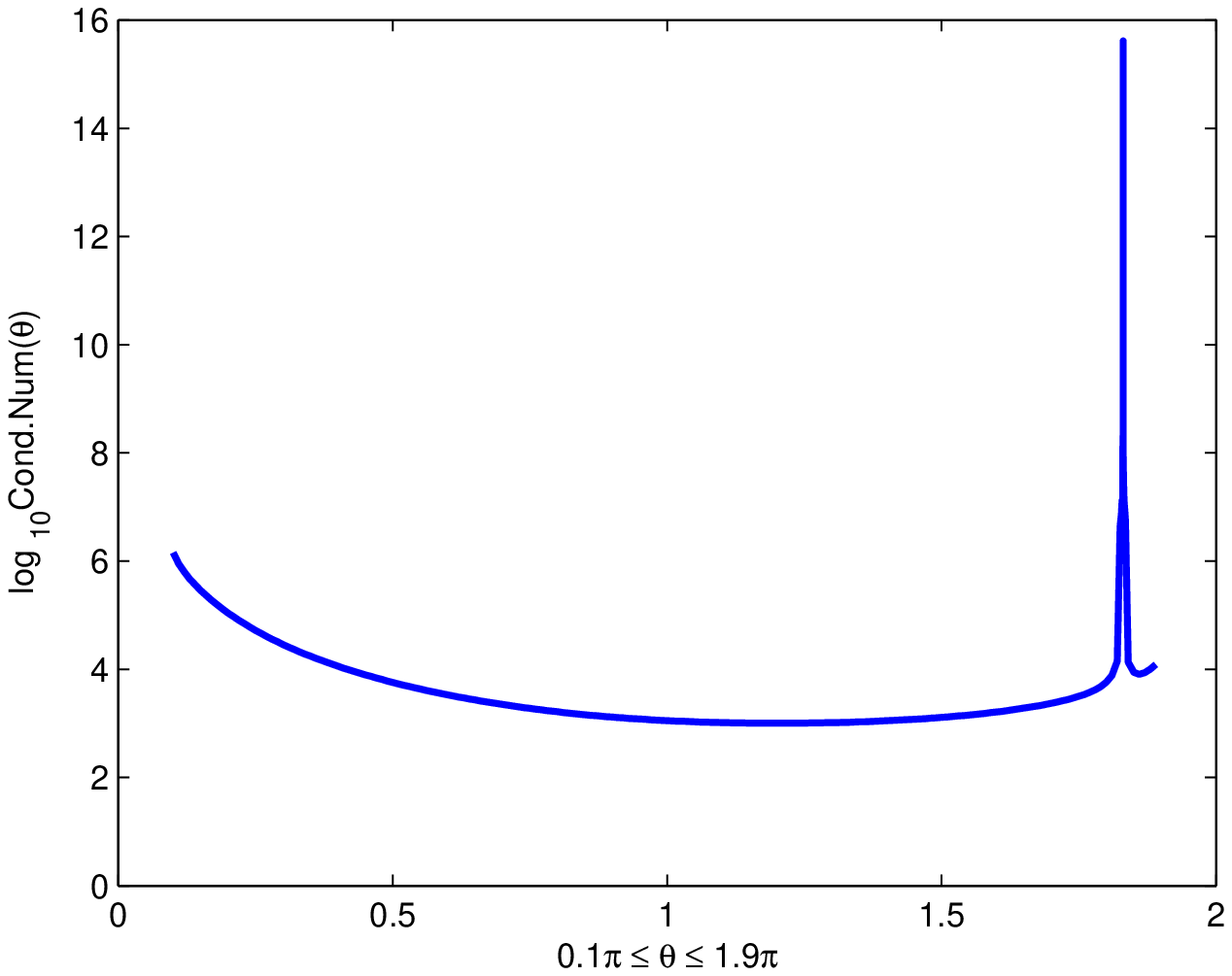}
\includegraphics[height=45mm,width=60mm]{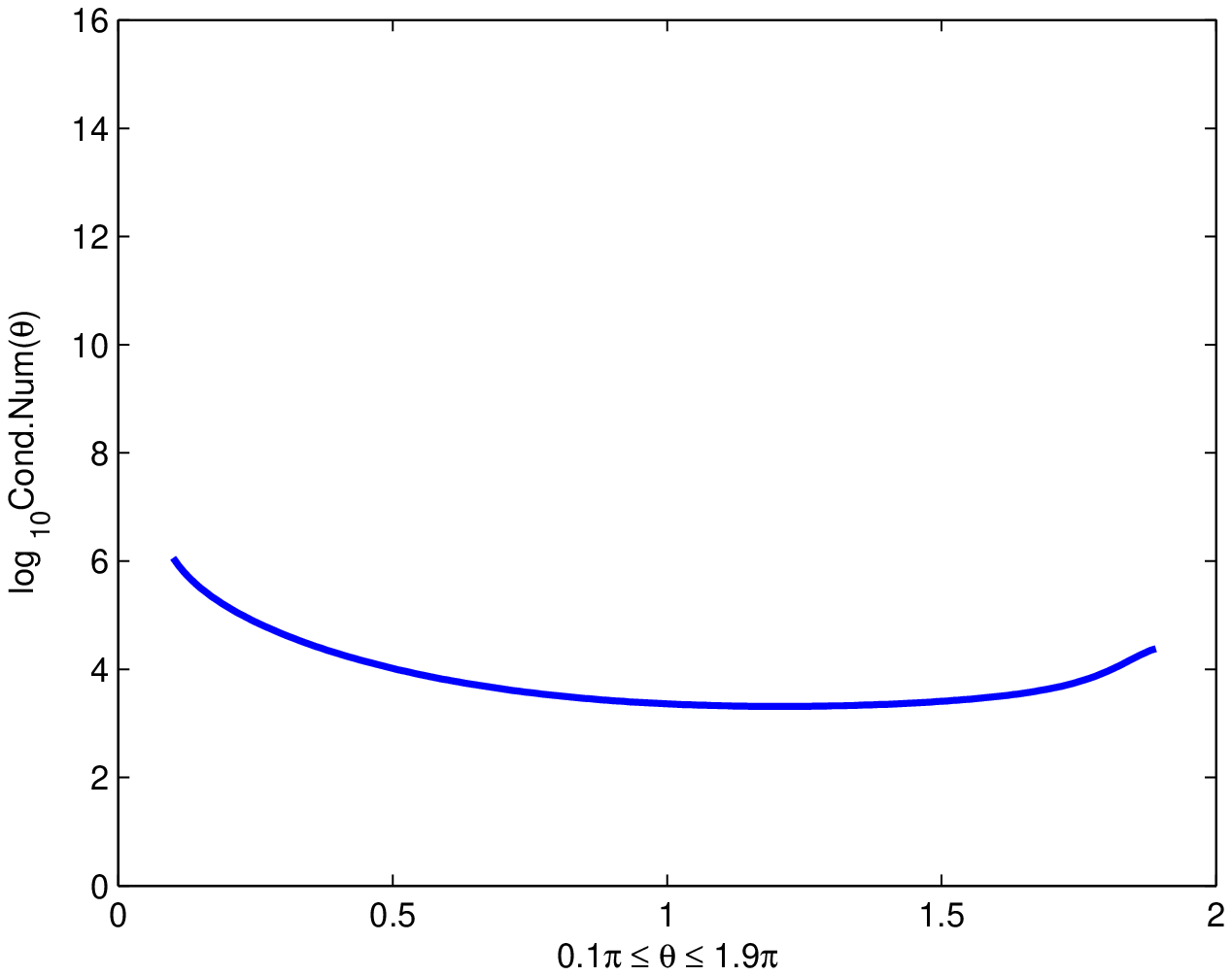}
\includegraphics[height=45mm,width=60mm]{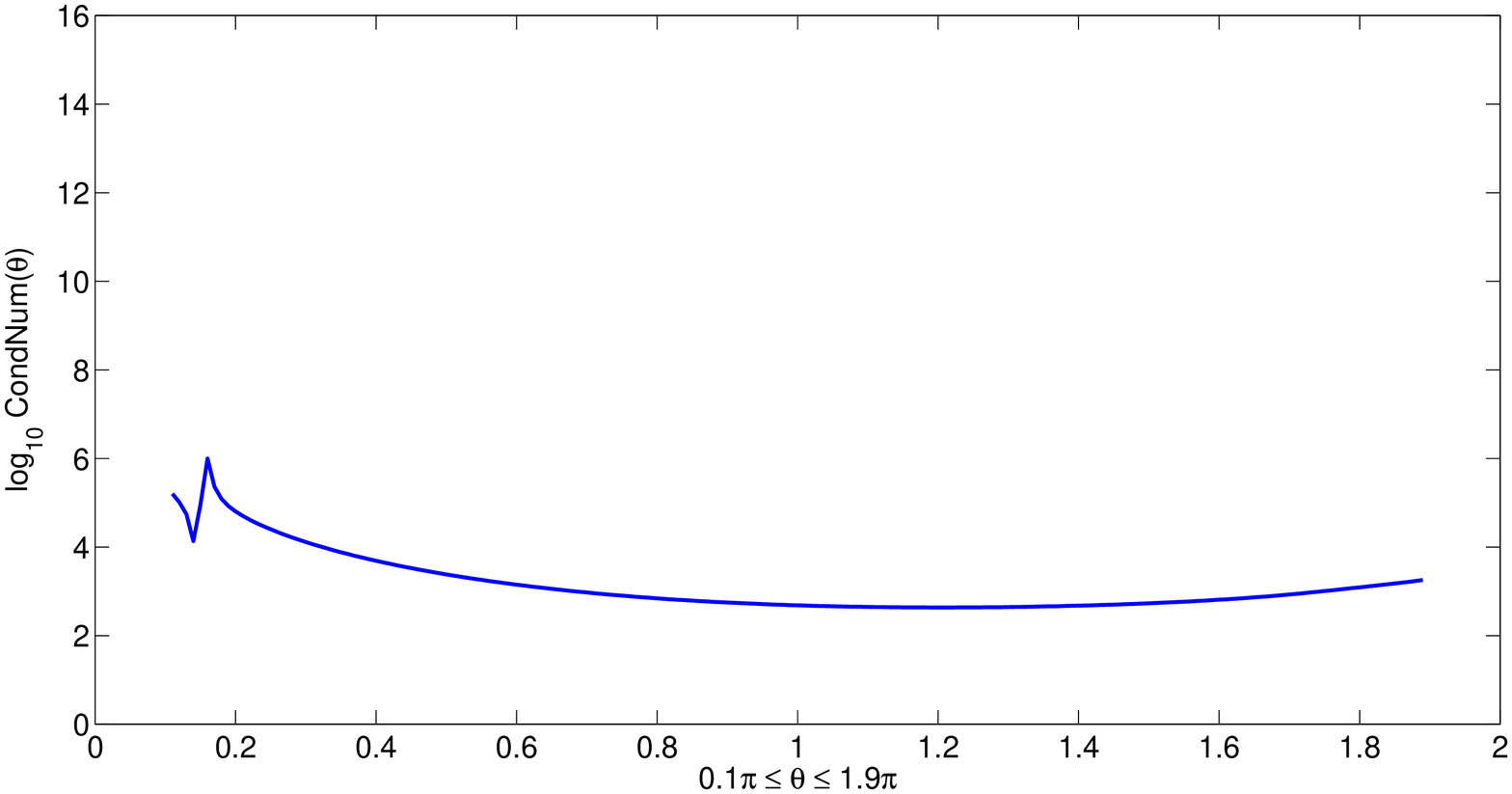}
\caption{\sf Condition numbers vs. opening angles in case $n=256$.
From row 1 to row 3: splines of degree 0, 1 and 2, respectively.
Left column: one-corner geometry, right column: two-corner
geometry.} \label{fig:256a}
\end{figure}

\begin{figure}[!ht]
\centering
\includegraphics[height=45mm,width=60mm]{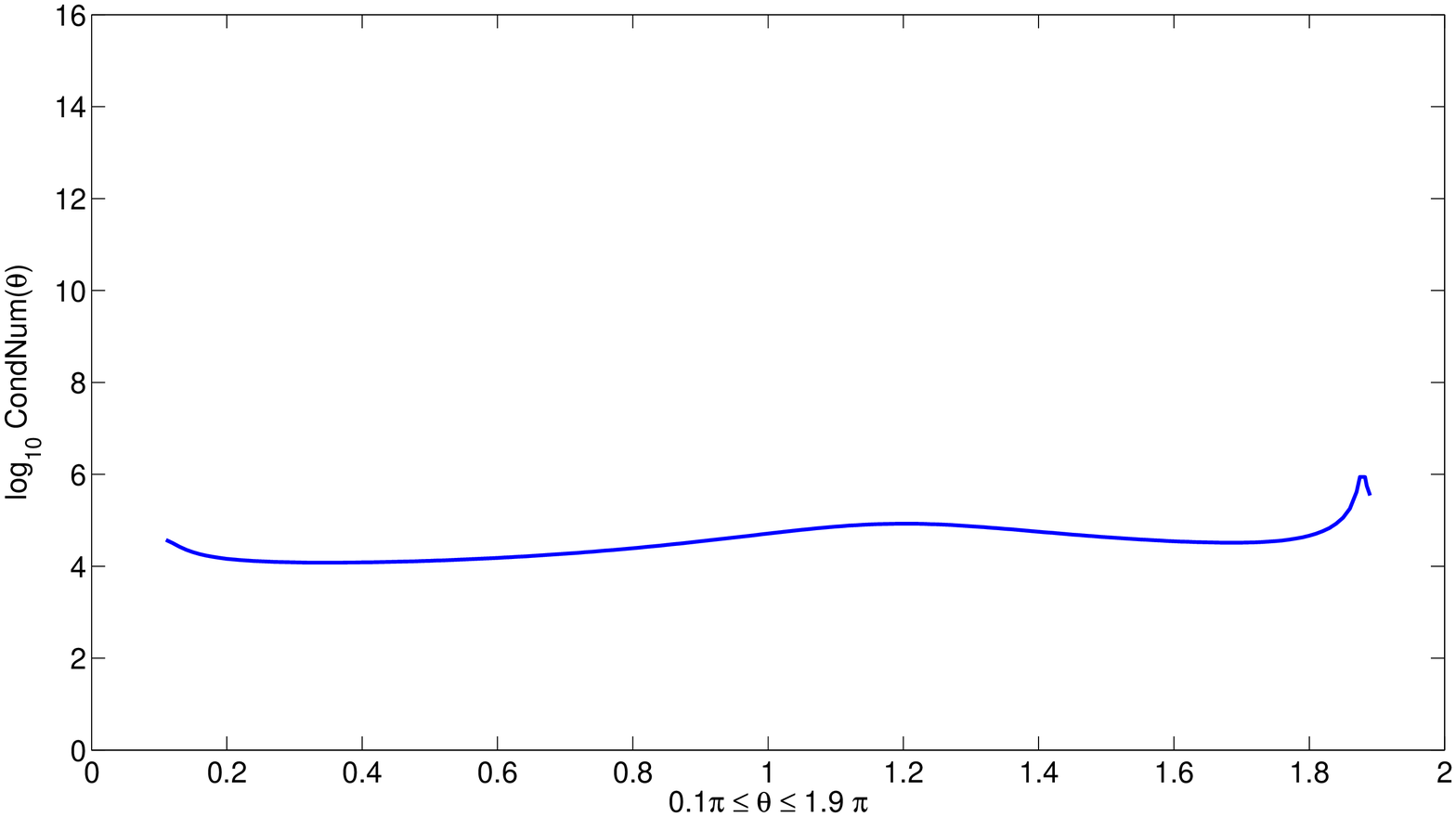}
\includegraphics[height=45mm,width=60mm]{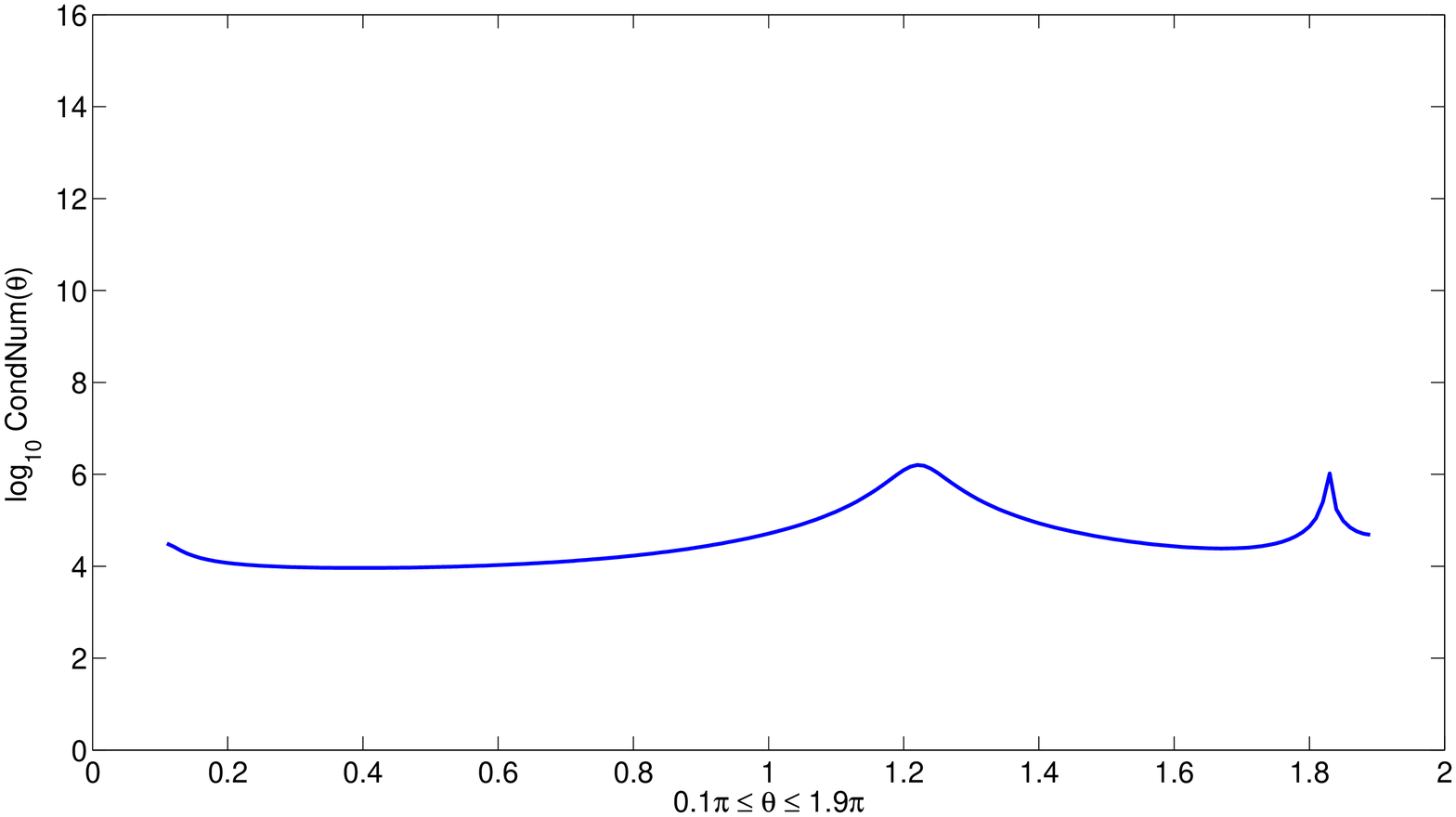}
\includegraphics[height=45mm,width=60mm]{glk7_d1_q1_256.eps}
\includegraphics[height=45mm,width=60mm]{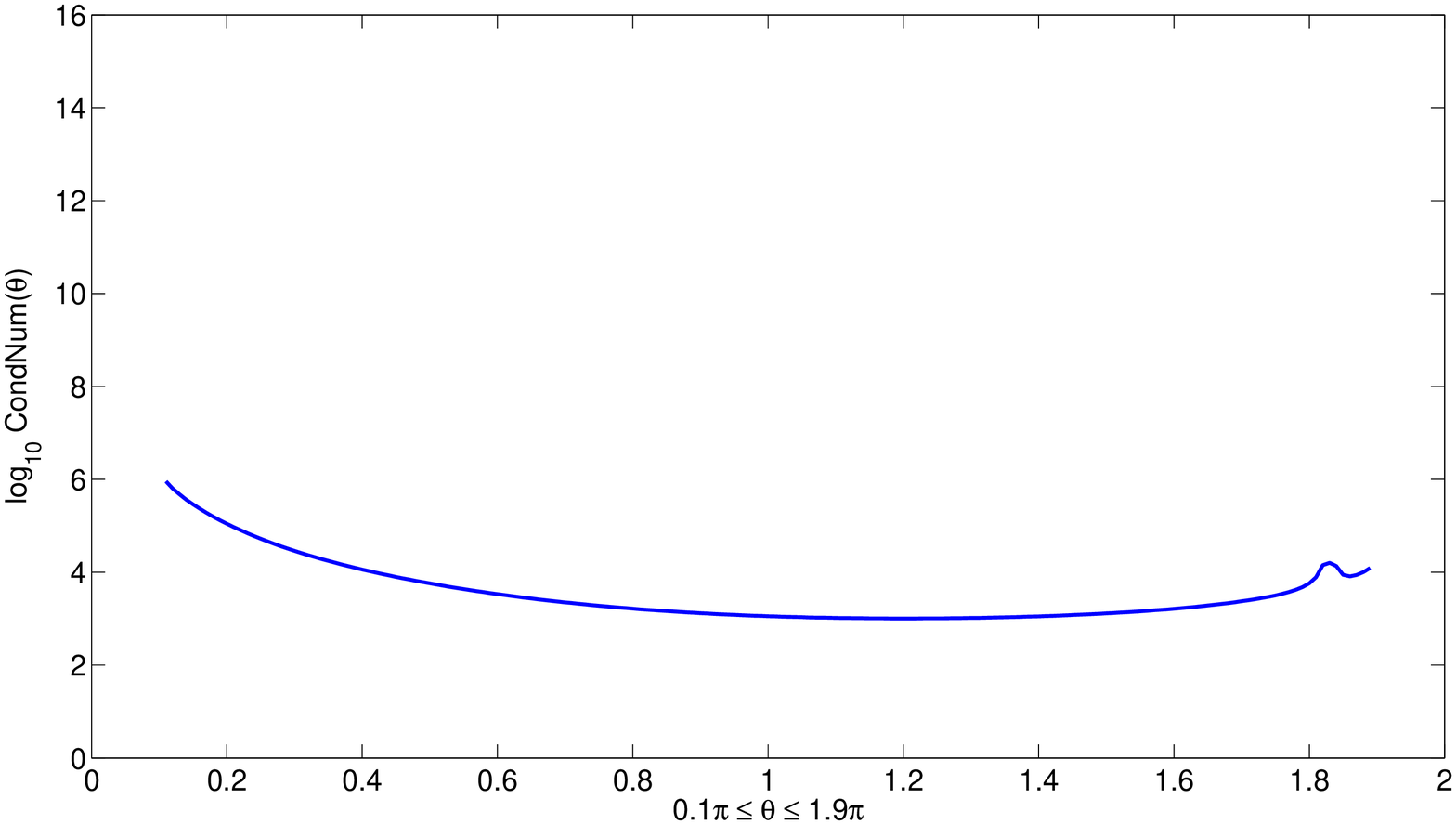}
\includegraphics[height=45mm,width=60mm]{glk7_d2_q1_256.eps}~%
\includegraphics[height=45mm,width=60mm]{glk7_d2_q2_256.eps}
\caption{\sf Condition numbers vs. opening angles in case $n=256$
and $n=512$ in neighbourhoods of suspicious points. From row 1 to
row 3: splines of degree 0, 1 and 2, respectively. Left column:
one-corner geometry, right column: two-corner geometry.}
\label{fig:256}
\end{figure}

Thus in order to determine the critical angles, i.e. the opening
angles $\theta$ for which the operators $R^\tau$ are not invertible,
one can consider the behaviour of the spline Galerkin methods on
special contours. A family of such contours $\Gamma_1^\theta$,
$\theta\in (0,2\pi)$,
$$
\Gamma_1^\theta:=\{ t\in\sC: t=\gamma _1(s) = \sin (\pi s) \exp (i
\theta (s-0.5)) ,\; s \in [0,1] \}
 $$
has been used in \cite{DH:2011b, DH:2013a} to study the local
operators of the Nystr\"om method for Sherman--Lauricella and
Muskhelishvili equations. Changing the parameter $\theta$ in the
interval $(0,2\pi)$, one obtains contours  located at the origin and
having only one corner of various magnitude. In the present paper,
we use the same contours to detect the critical angles of the spline
Galerkin methods. It is worth mentioning that the operator $R^\tau$
depends not only on $\theta_\tau$ but also on the angle $\beta_\tau$
between the right semi-tangent to the contour $\Gamma_1^\theta$ at
the point $\tau$ and the real line $\sR$. However, numerical
experiments conducted for both the Nystr\"om and spline Galerkin
methods show that, in fact, the angle $\beta_\tau$ does not
influence the invertibility of the operator $R^\tau$. This opens a
way for verifying the results obtained for contour $\Gamma_1^\theta$
by conducting similar tests for equations on contours with two or
more corner, all of the same magnitude. To this end, we will use
another contour $\Gamma_2^\theta$, which is the union of two
circular arcs with the parametrization
 \begin{align*}
\D \gamma_1(s)&=-0.5 \cot (0.5\theta )+0.5/\sin (0.5 \theta) \,
 \exp (i\theta (s-0.5)), \quad 0 \leq s \leq 1, \\[1ex]
\D \gamma_2(s)&=\phantom{-}0.5 \cot (0.5\theta )-0.5/\sin (0.5
\theta )\, \exp (i\theta (s-0.5)) , \quad  0 \leq s \leq 1.
\end{align*}
To find the angles of instability, the interval $[0.1\pi,\; 1.9\pi]$
has been divided by the points $\theta_k:=\pi(0.1+0.01k)$ and for
each opening angle $\theta_k$ we constructed the matrices of the
corresponding approximation operators for the Galerkin methods based
on the splines of degree $d=0, d=1$ and $d=2$.
\begin{figure}[!ht] \centering
\includegraphics[height=45mm]{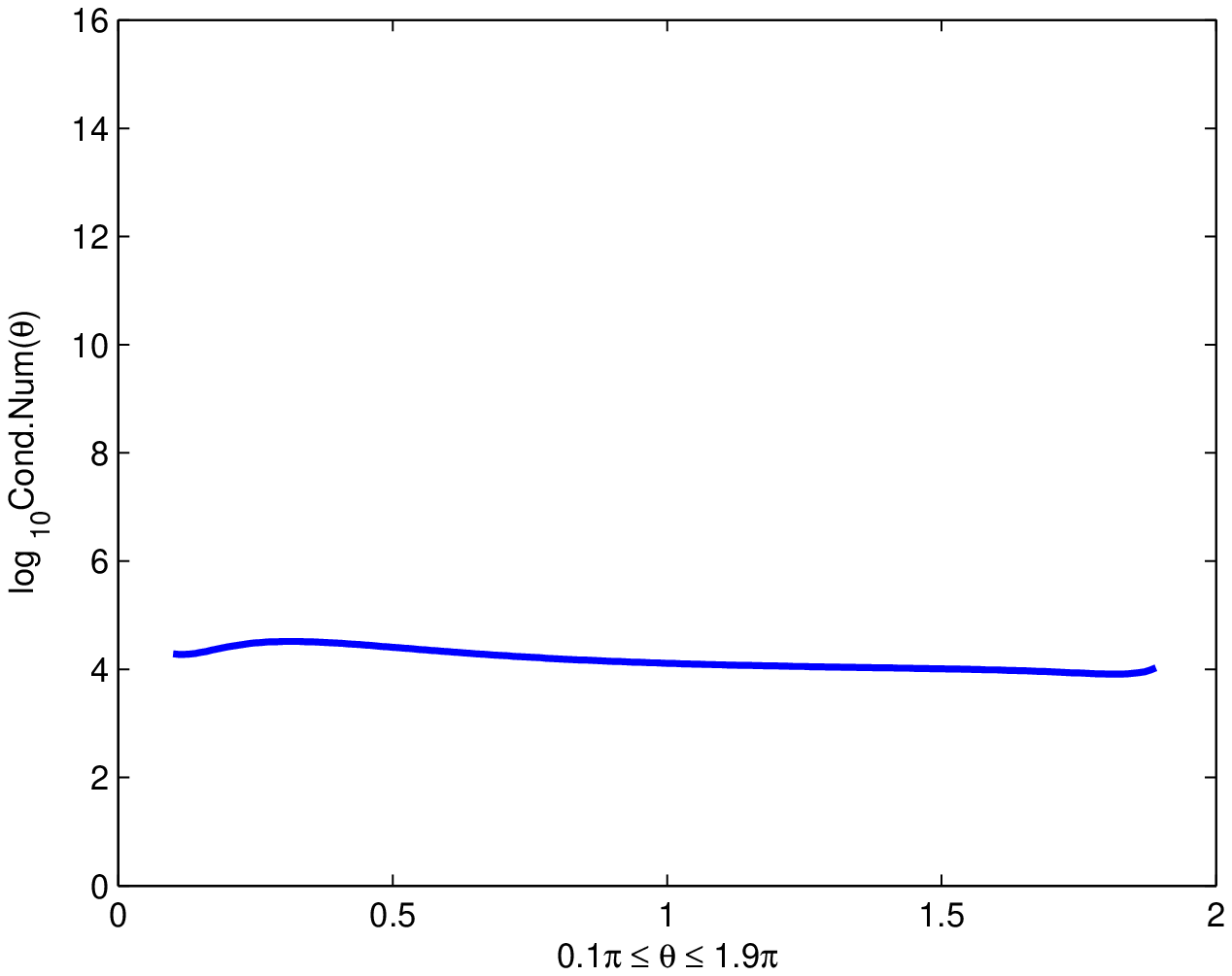}
 \includegraphics[height=45mm]{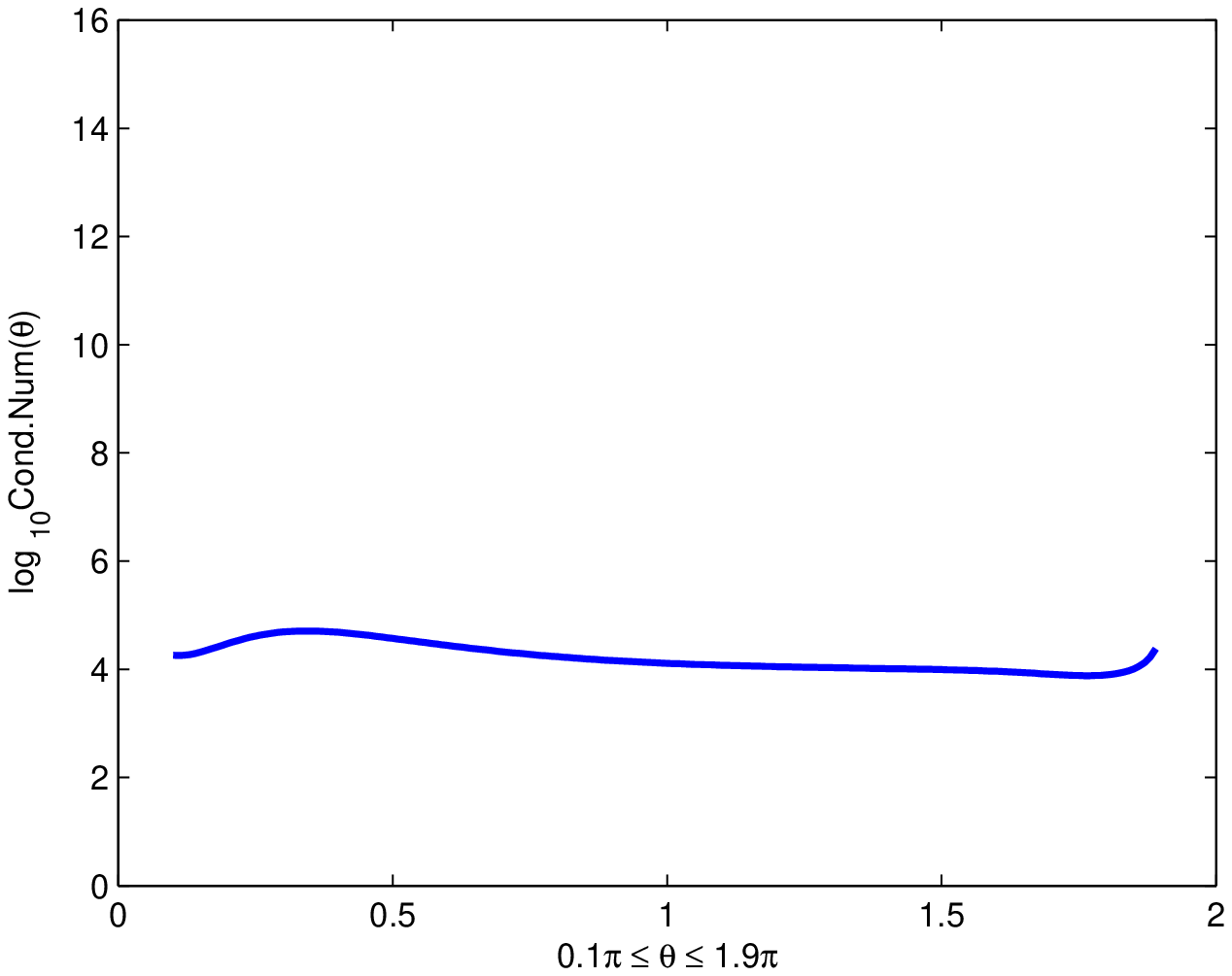}
 \includegraphics[height=45mm]{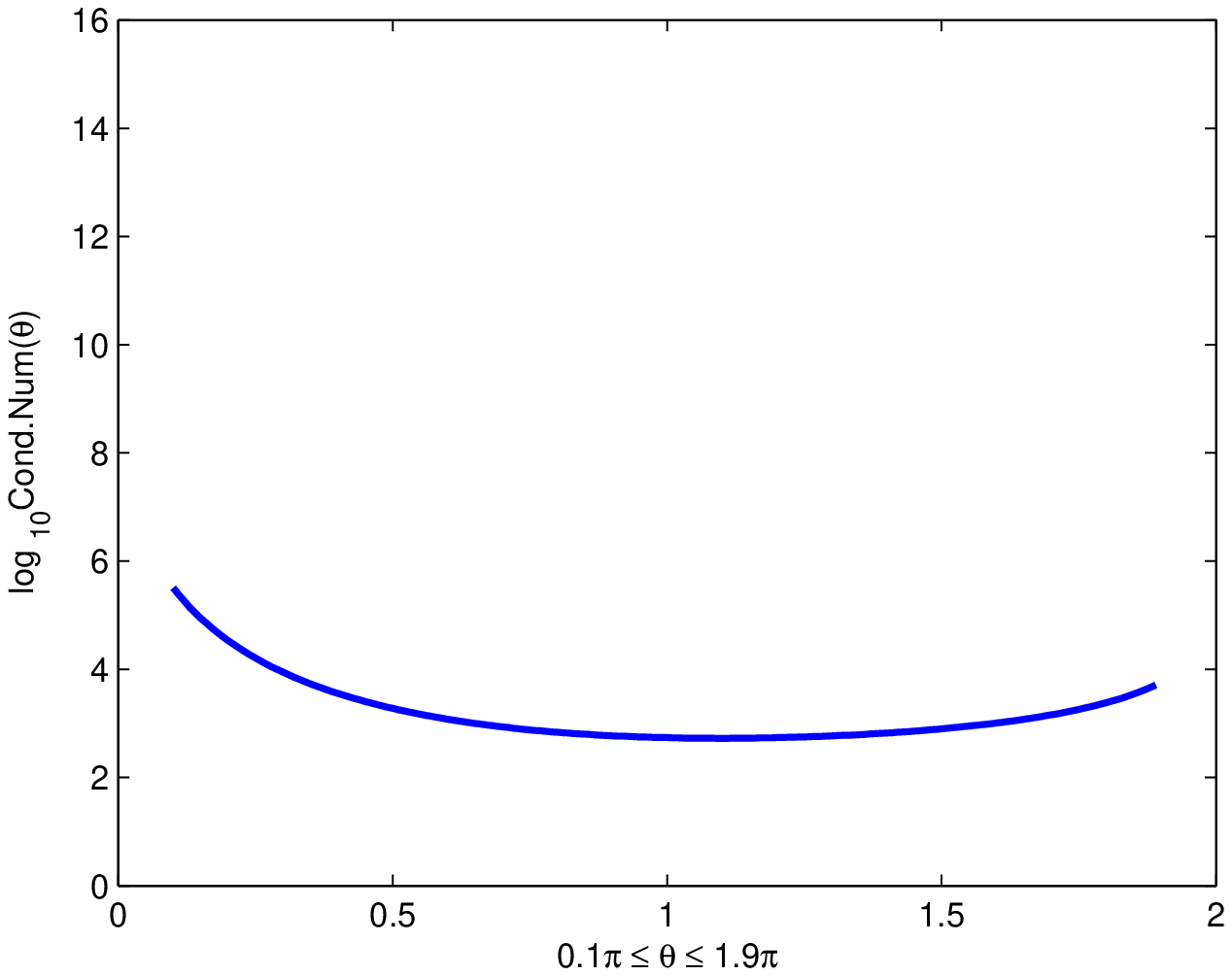}
  \includegraphics[height=45mm]{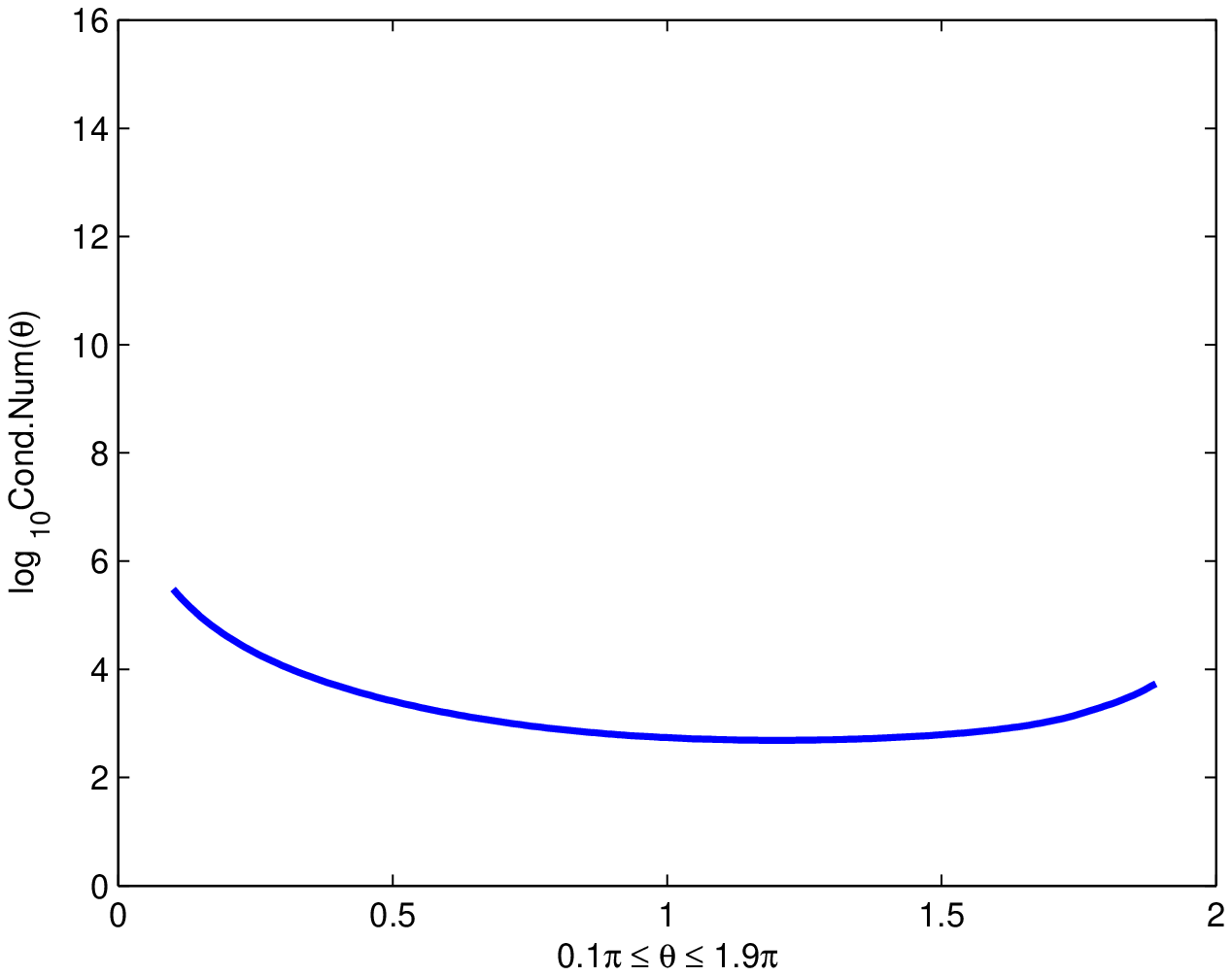}
 \includegraphics[height=45mm]{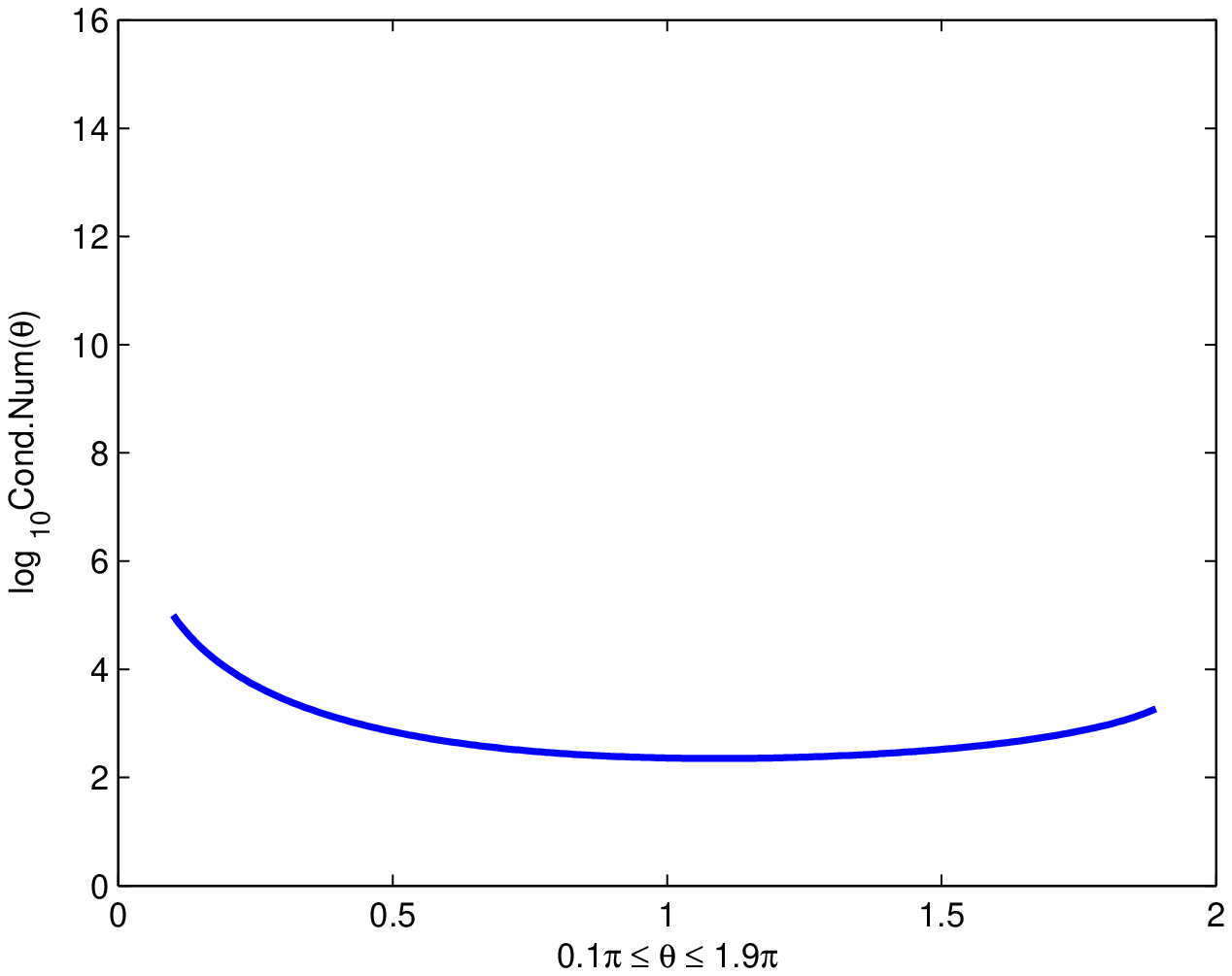}
  \includegraphics[height=45mm]{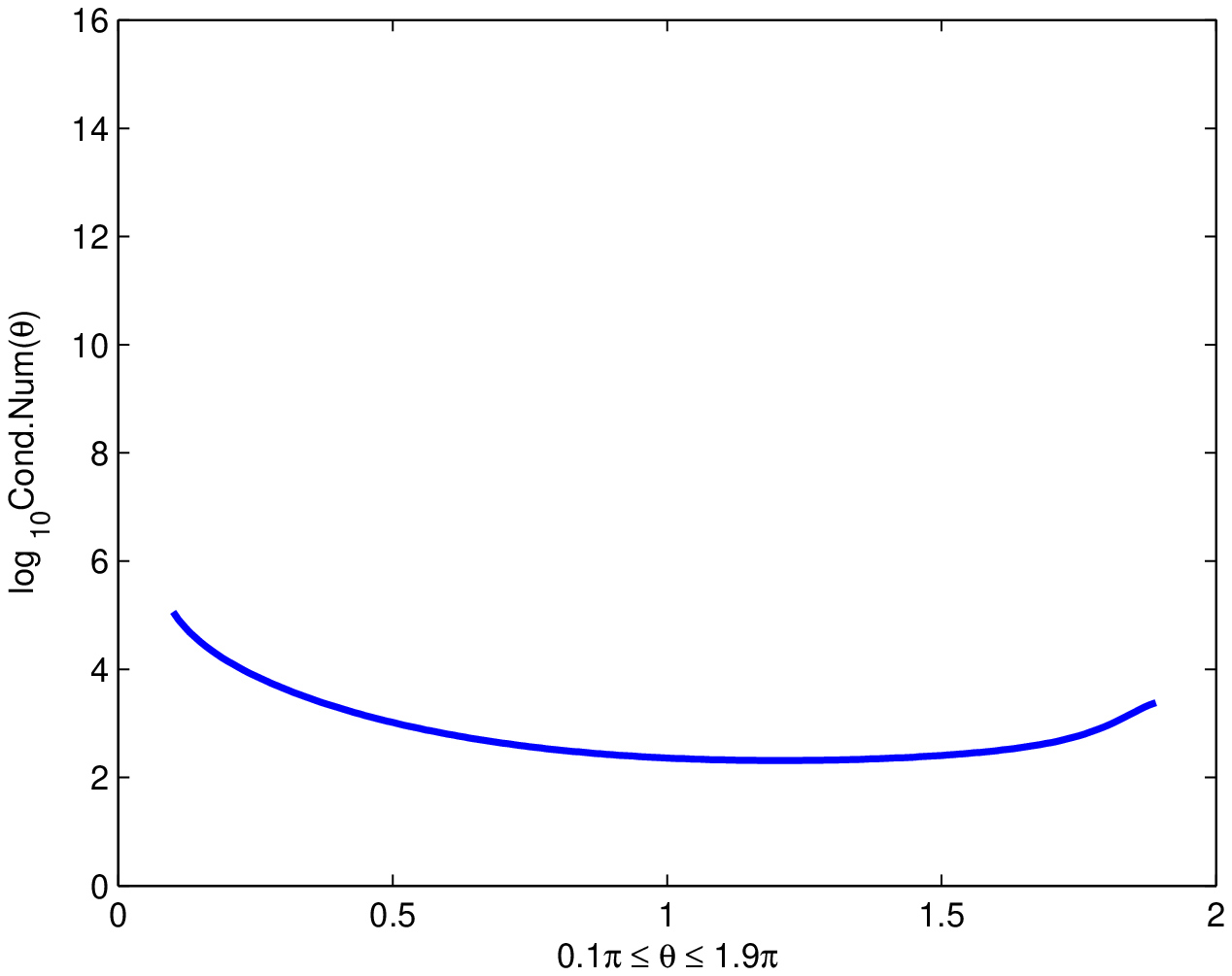}
   \caption{\sf Condition numbers vs. opening angles
   in case $n=128$. From row 1 to row 3: splines of degree 0, 1 and 2,
   respectively. Left column: one-corner geometry, right column: two-corner geometry.}
   \label{fig:128}
\end{figure}
Note that we consider Galerkin methods for two choices of $n$,
namely for $n=128$ and $n=256$, and the integrals arising in the
equation \eqref{eq12} and in the method \eqref{eq14} have been
approximated by quadrature formulas \eqref{quad1}, \eqref{quad2}.
Further, to verify the stability of the method, for each angle
$\theta_k$ we compute the condition numbers of the corresponding
matrices and the results of these computations are presented in
Figures \ref{fig:128}-\ref{fig:256}, where possible presence of
peaks might indicate critical angles. Thus it seems that inside of
the interval $(0.1\pi, 1.9\pi)$ neither of the Galerkin methods
based on splines of degree $0, 1$ or $2$ has critical angles. This
differs from the Nystr\"om method, where critical angles have been
discovered for both Sherman--Lauricella and Muskhelishvili equations
\cite{DH:2011b, DH:2013a}. Contrariwise, information about the
critical angles at the interval ends is not so conclusive. Thus in
the case $n=256$, the computation of the condition numbers for both
one and two corner geometry shows that for the Galerkin method based
on the splines of degree zero there can be a critical angle at the
right end of the interval mentioned.

For splines of the degree $d=0$ and $d=1$, the one and two corner
geometries give contradictory results (see Figure \ref{fig:256a}).
To clarify the situation one has to refine the mesh $\{\theta_k\}$
and essentially increase the dimension of the matrices used. Note
that while discovering a suspicious critical angle for $n=256$, we
refined the mesh $\{\theta_k\}$ in a neighbourhood of that angle by
reducing its step to $0.001 \pi$, and calculated the condition
numbers for the corresponding Galerkin methods with $n$ changed to
$512$. This allows us to show that, in fact, there are no critical
angles in the interval mentioned. However, the computing time
increases drastically.

The numerical experiments are performed in MATLAB environment
(version 7.9.0) and executed on an Acer Veriton M680 workstation
equipped with a Intel Core i7 vPro 870 Processor and 8GB of RAM, and
it took from one to two weeks of computer work in order to obtain
every single graph presented in Figure \ref{fig:128}, \ref{fig:256a}
or \ref{fig:256}.

 \end{document}